\documentclass[12pt,twoside]{amsart}
\usepackage[latin1]{inputenc}
\usepackage{amsmath, amsthm, amscd, amsfonts, amssymb, graphicx}
\usepackage[bookmarksnumbered, plainpages]{hyperref}

\newtheorem{thm}{Theorem}[section]
\newtheorem{cor}[thm]{Corollary}
\newtheorem{lem}[thm]{Lemma}

\newtheorem{defn}[thm]{Definition}

\numberwithin{equation}{section}

\begin{document}

\title{\bf The perturbation of the de Rham Hodge Operator and the Kastler-Kalau-Walze type theorem for manifolds with boundary}
\author{Siyao Liu  \hskip 0.4 true cm Tong Wu \hskip 0.4 true cm  Yong Wang$^{*}$}

\thanks{{\scriptsize
\hskip -0.4 true cm \textit{2010 Mathematics Subject Classification:}
53C40; 53C42.
\newline \textit{$^{*}$Corresponding author}}}

\maketitle

\begin{abstract}
 \indent In this paper, we give Lichnerowicz type formulas for the perturbation of the de Rham Hodge operator. We prove the Kastler-Kalau-Walze type theorems for the perturbation of the de Rham Hodge operator on 4-dimensional and 6-dimensional compact manifolds with or without boundary. Some concrete examples of the perturbation of the de Rham Hodge operator are provided for our main theorems.
\end{abstract}

\vskip 0.2 true cm


\pagestyle{myheadings}
\markboth{\rightline {\scriptsize  Liu}}
         {\leftline{\scriptsize }}

\bigskip
\bigskip


\section{ Introduction}
The noncommutative residue was found in \cite{Gu,Wo}. Since noncommutative residues are of great importance to the study of noncommutative geometry, more and more attention has been attached to the study of noncommutative residues.
By using the noncommutative residue, Connes derived a conformal 4-dimensional Polyakov action analogy \cite{Co1}. Connes showed us that the noncommutative residue on a compact manifold $M$ coincided with the Dixmiers trace on pseudodifferential operators of order $-{\rm {dim}}M$ in \cite{Co2}. Connes put forward that the noncommutative residue of the square of the inverse of the Dirac operator was proportioned to the Einstein-Hilbert action, which is called the Kastler-Kalau-Walze theorem now.
Kastler gave a brute-force proof of this theorem \cite{Ka}. In the same time, Kalau and Walze proved this theorem in the normal coordinates system \cite{KW}.
Ackermann proved this theorem by using the heat kernel expansion, in \cite{Ac}. The result of Connes was extended to the higher dimensional case \cite{U}.
Fedosov et al. gave the definition about the noncommutative residues on Boutet de Monvel algebra \cite{FGLS}.

On the other hand, Wang generalized the Connes' results to the case of manifolds with boundary in \cite{Wa1,Wa2}, and proved the Kastler-Kalau-Walze type theorem for the Dirac operator and the signature operator on lower-dimensional manifolds with boundary.
In \cite{Wa3,Wa4}, for the Dirac operator and the signature operator, Wang computed $\widetilde{{\rm Wres}}[\pi^+D^{-1}\circ\pi^+D^{-1}]$, in these cases the boundary term vanished. But for $\widetilde{{\rm Wres}}[\pi^+D^{-1}\circ\pi^+D^{-3}]$, authors got a nonvanishing boundary term \cite{Wa5}, and gave a theoretical explanation for gravitational action on boundary. In other words, Wang provided a kind of method to study the Kastler-Kalau-Walze type theorem for manifolds with boundary.
In \cite{WW}, Wei and Wang proved the Kastler-Kalau-Walze type theorem for modified Novikov operators on compact manifolds.
In \cite{WWW}, Wu-Wang-Wang obtained two Lichnerowicz type formulas for the Dirac-witten operators, and gave the proof of Kastler-Kalau-Walze type theorems for the Dirac-witten operators on 4-dimensional and 6-dimensional compact manifolds with (resp.without) boundary.
In \cite{Wa}, Wang proved a Kastler-Kalau-Walze type theorem for perturbations of the Dirac operators on compact manifolds with or without boundary.

The motivation of this paper is to prove the Kastler-Kalau-Walze type theorems for the perturbation of the de Rham Hodge operator. Specifically, we calculate $\widetilde{{\rm Wres}}[\pi^+{{D}_{A}}^{-1}\circ\pi^+({D}^*_{A})^{-1}],$ $\widetilde{{\rm Wres}}[\pi^+{{D}_{A}}^{-1}\circ\pi^+{{D}_{A}}^{-1}],$ $\widetilde{{\rm Wres}}[\pi^+{{D}_{A}}^{-1}\circ\pi^+({D}^*_{A}{D}_{A}{D}^*_{A})^{-1}]$ and $\widetilde{{\rm Wres}}[\pi^+{{D}_{A}}^{-1}\circ\pi^+{{D}_{A}}^{-3}],$ for ${D}_{A}$ see (2.5).

A brief description of the organization of this paper is as follows.
In Section 2, this paper will firstly introduce the basic notions of the perturbation of the de Rham Hodge operator, by means of which we can compute the Lichnerowicz formulas for the perturbation of the de Rham Hodge operator. And then we present the Kastler-Kalau-Walze type theorems  for the perturbation of the de Rham Hodge operator on $n$-dimensional compact manifolds without boundary.
In the next section, we calculate $\widetilde{{\rm Wres}}[\pi^+{{D}_{A}}^{-1}\circ\pi^+({D}^*_{A})^{-1}]$ and $\widetilde{{\rm Wres}}[\pi^+{{D}_{A}}^{-1}\circ\pi^+{{D}_{A}}^{-1}]$ for 4-dimensional compact manifolds with boundary.
In Section 4, we prove the Kastler-Kalau-Walze type theorems on 6-dimensional compact manifolds with boundary for the perturbation of the de Rham Hodge operator.

\vskip 1 true cm

\section{ The perturbation of the de Rham Hodge operator and its Lichnerowicz formulas}

We give some definitions and basic notions which we will use in this paper.
Let $M$ be a $n$-dimensional ($n\geq 3$) oriented compact Riemannian manifold with a Riemannian metric $g^{M}$. And let $\nabla^L$ be the Levi-Civita connection about $g^M$. In the local coordinates $\{x_i; 1\leq i\leq n\}$ and the
fixed orthonormal frame $\{\widetilde{e_1},\cdots,\widetilde{e_n}\}$, the connection matrix $(\omega_{s,t})$ is defined by
\begin{equation}
\nabla^L(\widetilde{e_1},\cdots,\widetilde{e_n})= (\widetilde{e_1},\cdots,\widetilde{e_n})(\omega_{s,t}).
\end{equation}
\indent Let $\epsilon (\widetilde{e_j*})$,~$\iota (\widetilde{e_j*})$ be the exterior and interior multiplications respectively, $\widetilde{e_j*}$ be the dual base of $\widetilde{e_j}$  and $c(\widetilde{e_j})$ be the Clifford action.
Suppose that $\partial_{i}$ is a natural local frame on $TM$
and $(g^{ij})_{1\leq i,j\leq n}$ is the inverse matrix associated to the metric
matrix  $(g_{ij})_{1\leq i,j\leq n}$ on $M$. Write
\begin{equation}
c(\widetilde{e_j})=\epsilon (\widetilde{e_j*})-\iota (\widetilde{e_j*});\
\overline{c}(\widetilde{e_j})=\epsilon (\widetilde{e_j*})+\iota (\widetilde{e_j*}),
\end{equation}
which satisfies
\begin{align}
&c(\widetilde{e_i})c(\widetilde{e_j})+c(\widetilde{e_j})c(\widetilde{e_i})=-2\delta_i^j;\\
&\overline{c}(\widetilde{e_i})c(\widetilde{e_j})+c(\widetilde{e_j})\overline{c}(\widetilde{e_i})=0;\nonumber\\
&\overline{c}(\widetilde{e_i})\overline{c}(\widetilde{e_j})+\overline{c}(\widetilde{e_j})\overline{c}(\widetilde{e_i})=2\delta_i^j.\nonumber
\end{align}
By \cite{Wa3}, we have the signature operator
\begin{align}
D&=d+\delta=\sum^n_{i=1}c(\widetilde{e_i})\bigg[\widetilde{e_i}+\frac{1}{4}\sum_{s,t}\omega_{s,t}
(\widetilde{e_i})[\overline{c}(\widetilde{e_s})\overline{c}(\widetilde{e_t})-c(\widetilde{e_s})c(\widetilde{e_t})]\bigg].
\end{align}

\indent  We define the perturbation of the de Rham Hodge operator ${D}_{A}$ as follows:
\begin{align}
{D}_{A}&=d+\delta+A\\
&=\sum^n_{i=1}c(\widetilde{e_i})\bigg[\widetilde{e}_i+\frac{1}{4}\sum_{s,t}\omega_{s,t}
(\widetilde{e_i})[\overline{c}(\widetilde{e_s})\overline{c}(\widetilde{e_t})-c(\widetilde{e_s})c(\widetilde{e_t})]\bigg]+A,\notag
\end{align}
then we have
\begin{align}
{D}^*_{A}&=d+\delta+A^*\\
&=\sum^n_{i=1}c(\widetilde{e_i})\bigg[\widetilde{e}_i+\frac{1}{4}\sum_{s,t}\omega_{s,t}
(\widetilde{e_i})[\overline{c}(\widetilde{e_s})\overline{c}(\widetilde{e_t})-c(\widetilde{e_s})c(\widetilde{e_t})]\bigg]+A^*,\notag
\end{align}
where
\begin{align}
A&=\overline{c}(X_1)\cdots\overline{c}(X_k)c(X_{k+1})\cdots c(X_{k+t});\\
A^*&=(-1)^tc(X_{k+t})\cdots c(X_{k+1})\overline{c}(X_k)\cdots\overline{c}(X_1),
\end{align}
where $X_1,\cdots,X_{k+t}$ are the smooth vector fields on $M$.\\
\indent By some simple calculations, we get Lichnerowicz formulas.
\begin{thm} The following equalities hold:
\begin{align}
{D}^2_{A}&=-[g^{ij}(\nabla_{\partial_{i}}^1\nabla_{\partial_{j}}^1-\nabla_{\nabla^{L}_{\partial_{i}}\partial_{j}}^1)]
-\frac{1}{8}\sum_{i,j,k,l}R_{ijkl}\overline{c}(\widetilde{e_i})\overline{c}(\widetilde{e_j})c(\widetilde{e_k})c(\widetilde{e_l})
+\frac{1}{4}s+A^2\\
&+\frac{1}{4}\sum_{i=1}^{n}[c(\widetilde{e_i})A+Ac(\widetilde{e_i})]^2
-\frac{1}{2}\sum_{j=1}^{n}[\nabla^{\bigwedge^*T^*M}_{\widetilde{e_j}}(A)c(\widetilde{e_j})-c(\widetilde{e_j})\nabla^{\bigwedge^*T^*M}_{\widetilde{e_j}}(A)],\nonumber\\
{D}^*_A{D}_{A}&=-[g^{ij}(\nabla_{\partial_{i}}^2\nabla_{\partial_{j}}^2-\nabla_{\nabla^{L}_{\partial_{i}}\partial_{j}}^2)]
-\frac{1}{8}\sum_{i,j,k,l}R_{ijkl}\overline{c}(\widetilde{e_i})\overline{c}(\widetilde{e_j})c(\widetilde{e_k})c(\widetilde{e_l})
+\frac{1}{4}s+A^*A\\
&+\frac{1}{4}\sum_{i=1}^{n}[c(\widetilde{e_i})A+A^*c(\widetilde{e_i})]^2
-\frac{1}{2}\sum_{j=1}^{n}[\nabla^{\bigwedge^*T^*M}_{\widetilde{e_j}}(A^*)c(\widetilde{e_j})-c(\widetilde{e_j})\nabla^{\bigwedge^*T^*M}_{\widetilde{e_j}}(A)], \notag
\end{align}
where $s$ is the scalar curvature.
\end{thm}
Now we prove the Theorem 2.1.
Let $M$ be a smooth compact oriented Riemannian $n$-dimensional manifolds without boundary and $N$ be a vector bundle on $M$.
We say that $P$ is a differential operator of Laplace type, if it has locally the form
\begin{equation}
P=-(g^{ij}\partial_i\partial_j+A^i\partial_i+B),
\end{equation}
where $\partial_{i}$  is a natural local frame on $TM$
and $(g^{ij})_{1\leq i,j\leq n}$ is the inverse matrix associated to the metric
matrix  $(g_{ij})_{1\leq i,j\leq n}$ on $M$,
 and $A^{i}$ and $B$ are smooth sections
of $\textrm{End}(N)$ on $M$ (endomorphism).
If $P$ satisfies the form (2.11), then there is a unique
connection $\nabla$ on $N$ and a unique endomorphism $E$ such that
 \begin{equation}
P=-[g^{ij}(\nabla_{\partial_{i}}\nabla_{\partial_{j}}- \nabla_{\nabla^{L}_{\partial_{i}}\partial_{j}})+E],
\end{equation}
where $\nabla^{L}$ is the Levi-Civita connection on $M$. Moreover
(with local frames of $T^{*}M$ and $N$), $\nabla_{\partial_{i}}=\partial_{i}+\omega_{i} $
and $E$ are related to $g^{ij}$, $A^{i}$ and $B$ through
 \begin{eqnarray}
&&\omega_{i}=\frac{1}{2}g_{ij}\big(A^{i}+g^{kl}\Gamma_{ kl}^{j} \texttt{id}\big),\\
&&E=B-g^{ij}\big(\partial_{i}(\omega_{j})+\omega_{i}\omega_{j}-\omega_{k}\Gamma_{ ij}^{k} \big),
\end{eqnarray}
where $\Gamma_{ kl}^{j}$ is the  Christoffel coefficient of $\nabla^{L}$.\\
\indent Let $g^{ij}=g(dx_{i},dx_{j})$, $\xi=\sum_{k}\xi_{j}dx_{j}$ and $\nabla^L_{\partial_{i}}\partial_{j}=\sum_{k}\Gamma_{ij}^{k}\partial_{k}$,  we denote that
\begin{align}
&\sigma_{i}=-\frac{1}{4}\sum_{s,t}\omega_{s,t}(\widetilde{e_i})c(\widetilde{e_s})c(\widetilde{e_t});\ a_{i}=\frac{1}{4}\sum_{s,t}\omega_{s,t}(\widetilde{e_i})\overline{c}(\widetilde{e_s})\overline{c}(\widetilde{e_t});\\
&\xi^{j}=g^{ij}\xi_{i};\ \Gamma^{k}=g^{ij}\Gamma_{ij}^{k};\ \sigma^{j}=g^{ij}\sigma_{i};\ a^{j}=g^{ij}a_{i}.\nonumber
\end{align}
\indent Then, the perturbation of the de Rham Hodge operator $D_{A}$ and ${D}^*_{A}$ can be written as
\begin{align}
{D}_{A}&=\sum^n_{i=1}c(\widetilde{e_i})[\widetilde{e_i}+\sigma_{i}+a_{i}]+A;\\
{D}^*_{A}&=\sum^n_{i=1}c(\widetilde{e_i})[\widetilde{e_i}+\sigma_{i}+a_{i}]+A^*.
\end{align}
\indent By \cite{Y}, the local expression of $(d+\delta)^2$ is
\begin{align}
(d+\delta)^2=-\Delta_0-\frac{1}{8}\sum_{i,j,k,l}R_{ijkl}\overline{c}(\widetilde{e_i})\overline{c}(\widetilde{e_j})c(\widetilde{e_k})c(\widetilde{e_l})+\frac{1}{4}s.
\end{align}
\indent By \cite{Ac}and\cite{Y}, we have
\begin{align}
-\Delta_0=-g^{ij}(\nabla_i^L\nabla_j^L-\Gamma_{ij}^k\nabla_k^L).
\end{align}
\indent We note that
\begin{align}
{D}^2_{A}&=(d+\delta)^2+(d+\delta)A+A(d+\delta)+A^2,
\end{align}
\begin{align}
(d+\delta)A+A(d+\delta)&=\sum_{i,j}g^{ij}[c(\partial_{i})A+Ac(\partial_{i})]\partial_{j}+\sum_{i,j}g^{ij}[c(\partial_{i})\partial_{j}(A)\\
&+c(\partial_{i})a_{j}A+c(\partial_{i})\sigma_{j}A+Ac(\partial_{i})a_{j}+Ac(\partial_{i})\sigma_{j}],\nonumber
\end{align}
then we obtain
\begin{align}
{D}^2_{A}&=-\sum_{i,j}g^{ij}[\partial_{i}\partial_{j}+2\sigma_{i}\partial_{j}+2a_{i}\partial_{j}-\Gamma_{ij}^{k}\partial_{k}
+(\partial_{i}\sigma_{j})+(\partial_{i}a_{j})+\sigma_{i}\sigma_{j}+\sigma_{i}a_{j}+a_{i}\sigma_{j}\\
&+a_{i}a_{j}-\Gamma_{ij}^{k}\sigma_{k}-\Gamma_{ij}^{k}a_{k}]+\sum_{i,j}g^{ij}[c(\partial_{i})\partial_{j}(A)+c(\partial_{i})a_{j}A+c(\partial_{i})\sigma_{j}A+Ac(\partial_{i})a_{j}\nonumber\\
&+Ac(\partial_{i})\sigma_{j}]+\sum_{i,j}g^{ij}[c(\partial_{i})A+Ac(\partial_{i})]\partial_{j}-\frac{1}{8}\sum_{i,j,k,l}R_{ijkl}\overline{c}(\widetilde{e_i})\overline{c}(\widetilde{e_j})c(\widetilde{e_k})c(\widetilde{e_l})\nonumber\\
&+\frac{1}{4}s+A^2.\nonumber
\end{align}
\indent Similarly, we have
\begin{align}
{D}^*_{A}{D}_{A}&=-\sum_{i,j}g^{ij}[\partial_{i}\partial_{j}+2\sigma_{i}\partial_{j}+2a_{i}\partial_{j}-\Gamma_{ij}^{k}\partial_{k}
+(\partial_{i}\sigma_{j})+(\partial_{i}a_{j})+\sigma_{i}\sigma_{j}+\sigma_{i}a_{j}+a_{i}\sigma_{j}\\
&+a_{i}a_{j}-\Gamma_{ij}^{k}\sigma_{k}-\Gamma_{ij}^{k}a_{k}]+\sum_{i,j}g^{ij}[c(\partial_{i})\partial_{j}(A)+c(\partial_{i})a_{j}A+c(\partial_{i})\sigma_{j}A+A^*c(\partial_{i})a_{j}\nonumber\\
&+A^*c(\partial_{i})\sigma_{j}]+\sum_{i,j}g^{ij}[c(\partial_{i})A+A^*c(\partial_{i})]\partial_{j}-\frac{1}{8}\sum_{i,j,k,l}R_{ijkl}\overline{c}(\widetilde{e_i})\overline{c}(\widetilde{e_j})c(\widetilde{e_k})c(\widetilde{e_l})\nonumber\\
&+\frac{1}{4}s+A^*A.\nonumber
\end{align}
\indent By (2.13), (2.14) and (2.22), we have
\begin{align}
(\omega_{i})_{{D}^2_{A}}&=\sigma_{i}+a_{i}-\frac{1}{2}[c(\partial_{i})A+Ac(\partial_{i})],
\end{align}
\begin{align}
E_{{D}^2_{A}}&=\frac{1}{8}\sum_{i,j,k,l}R_{ijkl}\overline{c}(\widetilde{e_i})\overline{c}(\widetilde{e_j})c(\widetilde{e_k})c(\widetilde{e_l})-\sum_{i=1}^{n}c(\partial_{i})\partial^{i}(A)-\sum_{i=1}^{n}c(\partial_{i})a^{i}A-\sum_{i=1}^{n}c(\partial_{i})\sigma^{i}A\\
&+\frac{1}{2}\sum_{j=1}^{n}\partial^{j}[c(\partial_{j})A+Ac(\partial_{j})]-\frac{1}{2}\sum_{k=1}^{n}\Gamma^k[c(\partial_{k})A+Ac(\partial_{k})]-\frac{1}{4}s-A^2\nonumber\\
&+\sum_{i,j}\frac{g^{ij}}{2}[c(\partial_{i})A+Ac(\partial_{i})](\sigma_{j}+a_{j})+\sum_{i,j}\frac{g^{ij}}{2}(\sigma_{i}+a_{i})[c(\partial_{j})A+Ac(\partial_{j})]\nonumber\\
&-\sum_{i,j}\frac{g^{ij}}{4}[c(\partial_{i})A+Ac(\partial_{i})][c(\partial_{j})A+Ac(\partial_{j})]-\sum_{i=1}^{n}Ac(\partial_{i})a^{i}-\sum_{i=1}^{n}Ac(\partial_{i})\sigma^{i}.\nonumber
\end{align}
\indent Since $E$ is globally
defined on $M$, taking normal coordinates at $x_0$, we have $\sigma^{i}(x_0)=0$, $a^{i}(x_0)=0$, $\partial^{j}[c(\partial_{j})](x_0)=0$, $\Gamma^k(x_0)=0$, $g^{ij}(x_0)=\delta^j_i$, so that
\begin{align}
E_{{D}^2_{A}}(x_0)&=\frac{1}{8}\sum_{i,j,k,l}R_{ijkl}\overline{c}(\widetilde{e_i})\overline{c}(\widetilde{e_j})c(\widetilde{e_k})c(\widetilde{e_l})-\frac{1}{4}\sum_{i=1}^{n}[c(\partial_{i})A+Ac(\partial_{i})]^2\\
&+\frac{1}{2}\sum_{j=1}^{n}[\partial^{j}(A)c(\partial_{j})-c(\partial_{j})\partial^{j}(A)]-\frac{1}{4}s-A^2\nonumber\\
&=\frac{1}{8}\sum_{i,j,k,l}R_{ijkl}\overline{c}(\widetilde{e_i})\overline{c}(\widetilde{e_j})c(\widetilde{e_k})c(\widetilde{e_l})-\frac{1}{4}\sum_{i=1}^{n}[c(\widetilde{e_i})A+Ac(\widetilde{e_i})]^2\nonumber\\
&+\frac{1}{2}\sum_{j=1}^{n}[\widetilde{e_j}(A)c(\widetilde{e_j})-c(\widetilde{e_j})\widetilde{e_j}(A)]-\frac{1}{4}s-A^2\nonumber\\
&=\frac{1}{8}\sum_{i,j,k,l}R_{ijkl}\overline{c}(\widetilde{e_i})\overline{c}(\widetilde{e_j})c(\widetilde{e_k})c(\widetilde{e_l})-\frac{1}{4}\sum_{i=1}^{n}[c(\widetilde{e_i})A+Ac(\widetilde{e_i})]^2\nonumber\\
&+\frac{1}{2}\sum_{j=1}^{n}[\nabla^{\bigwedge^*T^*M}_{\widetilde{e_j}}(A)c(\widetilde{e_j})-c(\widetilde{e_j})\nabla^{\bigwedge^*T^*M}_{\widetilde{e_j}}(A)]-\frac{1}{4}s-A^2.\nonumber
\end{align}
\indent Similarly, we have
\begin{align}
E_{{D}^*_{A}{D}_{A}}(x_0)&=\frac{1}{8}\sum_{i,j,k,l}R_{ijkl}\overline{c}(\widetilde{e_i})\overline{c}(\widetilde{e_j})c(\widetilde{e_k})c(\widetilde{e_l})-\frac{1}{4}\sum_{i=1}^{n}[c(\widetilde{e_i})A+A^*c(\widetilde{e_i})]^2\\
&+\frac{1}{2}\sum_{j=1}^{n}[\nabla^{\bigwedge^*T^*M}_{\widetilde{e_j}}(A^*)c(\widetilde{e_j})-c(\widetilde{e_j})\nabla^{\bigwedge^*T^*M}_{\widetilde{e_j}}(A)]-\frac{1}{4}s-A^*A.\nonumber
\end{align}
By (2.12), we get Theorem 2.1.\\
\indent According to the detailed descriptions in \cite{Ac}, we know that the noncommutative residue of a generalized laplacian $\widetilde{\Delta}$ is expressed as
\begin{equation}
(n-2)\Phi_{2}(\widetilde{\Delta})=(4\pi)^{-\frac{n}{2}}\Gamma(\frac{n}{2})\widetilde{res}(\widetilde{\Delta}^{-\frac{n}{2}+1}),
\end{equation}
where $\Phi_{2}(\widetilde{\Delta})$ denotes the integral over the diagonal part of the second
coefficient of the heat kernel expansion of $\widetilde{\Delta}$.
Now let $\widetilde{\Delta}={D}^2_{A}$. Since ${D}^2_{A}$ is a generalized laplacian, we can suppose ${D}^2_{A}=\Delta-E$, then, we have
\begin{align}
{\rm Wres}({D}^2_{A})^{-\frac{n-2}{2}}
=\frac{(n-2)(4\pi)^{\frac{n}{2}}}{(\frac{n}{2}-1)!}\int_{M}{\rm tr}(\frac{1}{6}s+E_{{D}^2_{A}})d{\rm Vol_{M} },
\end{align}
\begin{align}
{\rm Wres}({D}^*_{A}{D}_{A})^{-\frac{n-2}{2}}
=\frac{(n-2)(4\pi)^{\frac{n}{2}}}{(\frac{n}{2}-1)!}\int_{M}{\rm tr}(\frac{1}{6}s+E_{{D}^*_{A}{D}_{A}})d{\rm Vol_{M} },
\end{align}
where ${\rm Wres}$ denote the noncommutative residue.\\
By applying the formulae shown in (2.26), (2.27), (2.29) and (2.30), we get:
\begin{thm} If $M$ is a $n$-dimensional compact oriented manifolds without boundary, we have the following:
\begin{align}
&{\rm Wres}({D}^2_{A})^{-\frac{n-2}{2}}
=\frac{(n-2)(4\pi)^{\frac{n}{2}}}{(\frac{n}{2}-1)!}\int_{M}{\rm tr}\bigg(-\frac{1}{12}s+(\frac{n}{2}-1)A^2-\frac{1}{2}\sum_{j=1}^{n}Ac(\widetilde{e_j})Ac(\widetilde{e_j})\bigg)d{\rm Vol_{M}},\\
&{\rm Wres}({D}^*_{A}{D}_{A})^{-\frac{n-2}{2}}=\frac{(n-2)(4\pi)^{\frac{n}{2}}}{(\frac{n}{2}-1)!}\int_{M}{\rm tr}\bigg(-\frac{1}{12}s+(\frac{n}{2}-1)A^*A-\frac{1}{4}\sum_{j=1}^{n}Ac(\widetilde{e_j})Ac(\widetilde{e_j})\\
&-\frac{1}{4}\sum_{j=1}^{n}A^*c(\widetilde{e_j})A^*c(\widetilde{e_j})+\frac{1}{2}\sum_{j=1}^{n}\nabla^{\bigwedge^*T^*M}_{\widetilde{e_j}}(A^*)c(\widetilde{e_j})-\frac{1}{2}\sum_{j=1}^{n}c(\widetilde{e_j})\nabla^{\bigwedge^*T^*M}_{\widetilde{e_j}}(A)\bigg)d{\rm Vol_{M}},\nonumber
\end{align}
where $s$ is the scalar curvature.
\end{thm}

\vskip 1 true cm

\section{ A Kastler-Kalau-Walze type theorem for $4$-dimensional manifolds with boundary}
\indent Firstly, we explain the basic notions of  Boutet de Monvel's calculus and the definition of the noncommutative residue for manifolds with boundary that will be used throughout the paper. For the details, see Ref.\cite{Wa3}.\\
\indent Let $U\subset M$ be a collar neighborhood of $\partial M$ which is diffeomorphic with $\partial M\times [0,1)$. By the definition of $h(x_n)\in C^{\infty}([0,1))$
and $h(x_n)>0$, there exists $\widehat{h}\in C^{\infty}((-\varepsilon,1))$ such that $\widehat{h}|_{[0,1)}=h$ and $\widehat{h}>0$ for some
sufficiently small $\varepsilon>0$. Then there exists a metric $g'$ on $\widetilde{M}=M\bigcup_{\partial M}\partial M\times
(-\varepsilon,0]$ which has the form on $U\bigcup_{\partial M}\partial M\times (-\varepsilon,0 ]$
\begin{equation}
g'=\frac{1}{\widehat{h}(x_{n})}g^{\partial M}+dx _{n}^{2} ,
\end{equation}
such that $g'|_{M}=g$. We fix a metric $g'$ on the $\widetilde{M}$ such that $g'|_{M}=g$.

We define the Fourier transformation $F'$  by
\begin{equation}
F':L^2({\bf R}_t)\rightarrow L^2({\bf R}_v);~F'(u)(v)=\int e^{-ivt}u(t)dt\\
\end{equation}
and let
\begin{equation}
r^{+}:C^\infty ({\bf R})\rightarrow C^\infty (\widetilde{{\bf R}^+});~ f\rightarrow f|\widetilde{{\bf R}^+};~
\widetilde{{\bf R}^+}=\{x\geq0;x\in {\bf R}\}.
\end{equation}
 where $\Phi({\bf R})$
denotes the Schwartz space and $\Phi(\widetilde{{\bf R}^+}) =r^+\Phi({\bf R})$, $\Phi(\widetilde{{\bf R}^-}) =r^-\Phi({\bf R})$.\\
\indent We define $H^+=F'(\Phi(\widetilde{{\bf R}^+}));~ H^-_0=F'(\Phi(\widetilde{{\bf R}^-}))$ which satisfies
$H^+\bot H^-_0$. We have the following
 property: $h\in H^+~(H^-_0)$ if and only if $h\in C^\infty({\bf R})$ which has an analytic extension to the lower (upper) complex
half-plane $\{{\rm Im}\xi<0\}~(\{{\rm Im}\xi>0\})$ such that for all nonnegative integer $l$,
 \begin{equation}
\frac{d^{l}h}{d\xi^l}(\xi)\sim\sum^{\infty}_{k=1}\frac{d^l}{d\xi^l}(\frac{c_k}{\xi^k}),
\end{equation}
as $|\xi|\rightarrow +\infty,{\rm Im}\xi\leq0~({\rm Im}\xi\geq0)$.\\
\indent Let $H'$ be the space of all polynomials and $H^-=H^-_0\bigoplus H';~H=H^+\bigoplus H^-.$ Denote by $\pi^+~(\pi^-)$ respectively the projection on $H^+~(H^-)$. For calculations, we take $H=\widetilde H=\{$rational functions having no poles on the real axis$\}$ ($\tilde{H}$ is a dense set in the topology of $H$). Then on $\tilde{H}$,
 \begin{equation}
\pi^+h(\xi_0)=\frac{1}{2\pi i}\lim_{u\rightarrow 0^{-}}\int_{\Gamma^+}\frac{h(\xi)}{\xi_0+iu-\xi}d\xi,
\end{equation}
where $\Gamma^+$ is a Jordan close curve
included ${\rm Im}(\xi)>0$ surrounding all the singularities of $h$ in the upper half-plane and
$\xi_0\in {\bf R}$. Similarly, define $\pi'$ on $\tilde{H}$,
\begin{equation}
\pi'h=\frac{1}{2\pi}\int_{\Gamma^+}h(\xi)d\xi.
\end{equation}
So, $\pi'(H^-)=0$. For $h\in H\bigcap L^1({\bf R})$, $\pi'h=\frac{1}{2\pi}\int_{{\bf R}}h(v)dv$ and for $h\in H^+\bigcap L^1({\bf R})$, $\pi'h=0$.

Let $M$ be a $n$-dimensional compact oriented manifold with boundary $\partial M$.
Denote by $\mathcal{B}$ Boutet de Monvel's algebra, we recall the main theorem in \cite{FGLS,Wa3}.
\begin{thm}\label{th:32}{\rm\cite{FGLS}}{\bf(Fedosov-Golse-Leichtnam-Schrohe)}
 Let $X$ and $\partial X$ be connected, ${\rm dim}X=n\geq3$,
 $A=\left(\begin{array}{lcr}\pi^+P+G &   K \\
T &  S    \end{array}\right)$ $\in \mathcal{B}$ , and denote by $p$, $b$ and $s$ the local symbols of $P,G$ and $S$ respectively.
 Define:
 \begin{align}
{\rm{\widetilde{Wres}}}(A)&=\int_X\int_{\bf S}{\rm{tr}}_E\left[p_{-n}(x,\xi)\right]\sigma(\xi)dx \\
&+2\pi\int_ {\partial X}\int_{\bf S'}\left\{{\rm tr}_E\left[({\rm{tr}}b_{-n})(x',\xi')\right]+{\rm{tr}}
_F\left[s_{1-n}(x',\xi')\right]\right\}\sigma(\xi')dx',\nonumber
\end{align}
Then~~ a) ${\rm \widetilde{Wres}}([A,B])=0 $, for any
$A,B\in\mathcal{B}$;~~ b) It is a unique continuous trace on
$\mathcal{B}/\mathcal{B}^{-\infty}$.
\end{thm}

\begin{defn}{\rm\cite{Wa3} }
Lower dimensional volumes of spin manifolds with boundary are defined by
 \begin{equation}
{\rm Vol}^{(p_1,p_2)}_nM:= \widetilde{{\rm Wres}}[\pi^+D^{-p_1}\circ\pi^+D^{-p_2}],
\end{equation}
\end{defn}
 By \cite{Wa3}, we get
\begin{align}
\widetilde{{\rm Wres}}[\pi^+D^{-p_1}\circ\pi^+D^{-p_2}]=\int_M\int_{|\xi|=1}{\rm
trace}_{\wedge^*T^*M}[\sigma_{-n}(D^{-p_1-p_2})]\sigma(\xi)dx+\int_{\partial M}\Phi
\end{align}
and
\begin{align}
\Phi&=\int_{|\xi'|=1}\int^{+\infty}_{-\infty}\sum^{\infty}_{j, k=0}\sum\frac{(-i)^{|\alpha|+j+k+1}}{\alpha!(j+k+1)!}
\times {\rm trace}_{\wedge^*T^*M}[\partial^j_{x_n}\partial^\alpha_{\xi'}\partial^k_{\xi_n}\sigma^+_{r}(D^{-p_1})(x',0,\xi',\xi_n)
\\
&\times\partial^\alpha_{x'}\partial^{j+1}_{\xi_n}\partial^k_{x_n}\sigma_{l}(D^{-p_2})(x',0,\xi',\xi_n)]d\xi_n\sigma(\xi')dx',\nonumber
\end{align}
 where the sum is taken over $r+l-k-|\alpha|-j-1=-n,~~r\leq -p_1,l\leq -p_2$.

 Since $[\sigma_{-n}(D^{-p_1-p_2})]|_M$ has the same expression as $\sigma_{-n}(D^{-p_1-p_2})$ in the case of manifolds without
boundary, so locally we can compute the first term by \cite{Ka}, \cite{KW}, \cite{Wa3}, \cite{Po}.

For any fixed point $x_0\in\partial M$, we choose the normal coordinates
$U$ of $x_0$ in $\partial M$ (not in $M$) and compute $\Phi(x_0)$ in the coordinates $\widetilde{U}=U\times [0,1)\subset M$ and the
metric $\frac{1}{h(x_n)}g^{\partial M}+dx_n^2.$ The dual metric of $g^M$ on $\widetilde{U}$ is ${h(x_n)}g^{\partial M}+dx_n^2.$  Write
$g^M_{ij}=g^M(\frac{\partial}{\partial x_i},\frac{\partial}{\partial x_j});~ g_M^{ij}=g^M(dx_i,dx_j)$, then

\begin{equation}
[g^M_{i,j}]= \left[\begin{array}{lcr}
  \frac{1}{h(x_n)}[g_{i,j}^{\partial M}]  & 0  \\
   0  &  1
\end{array}\right];~~~
[g_M^{i,j}]= \left[\begin{array}{lcr}
  h(x_n)[g^{i,j}_{\partial M}]  & 0  \\
   0  &  1
\end{array}\right]
\end{equation}
and
\begin{equation}
\partial_{x_s}g_{ij}^{\partial M}(x_0)=0, 1\leq i,j\leq n-1; ~~~g_{ij}^M(x_0)=\delta_{ij}.
\end{equation}
\indent We review the following three lemmas.
\begin{lem}{\rm \cite{Wa3}}\label{le:32}
With the metric $g^{M}$ on $M$ near the boundary
\begin{eqnarray}
\partial_{x_j}(|\xi|_{g^M}^2)(x_0)&=&\left\{
       \begin{array}{c}
        0,  ~~~~~~~~~~ ~~~~~~~~~~ ~~~~~~~~~~~~~{\rm if }~j<n, \\[2pt]
       h'(0)|\xi'|^{2}_{g^{\partial M}},~~~~~~~~~~~~~~~~~~~~{\rm if }~j=n;
       \end{array}
    \right. \\
\partial_{x_j}[c(\xi)](x_0)&=&\left\{
       \begin{array}{c}
      0,  ~~~~~~~~~~ ~~~~~~~~~~ ~~~~~~~~~~~~~{\rm if }~j<n,\\[2pt]
\partial x_{n}(c(\xi'))(x_{0}), ~~~~~~~~~~~~~~~~~{\rm if }~j=n,
       \end{array}
    \right.
\end{eqnarray}
where $\xi=\xi'+\xi_{n}dx_{n}$.
\end{lem}
\begin{lem}{\rm \cite{Wa3}}\label{le:32}With the metric $g^{M}$ on $M$ near the boundary
\begin{align}
\omega_{s,t}(\widetilde{e_i})(x_0)&=\left\{
       \begin{array}{c}
        \omega_{n,i}(\widetilde{e_i})(x_0)=\frac{1}{2}h'(0),  ~~~~~~~~~~ ~~~~~~~~~~~{\rm if }~s=n,t=i,i<n; \\[2pt]
       \omega_{i,n}(\widetilde{e_i})(x_0)=-\frac{1}{2}h'(0),~~~~~~~~~~~~~~~~~~~{\rm if }~s=i,t=n,i<n;\\[2pt]
    \omega_{s,t}(\widetilde{e_i})(x_0)=0,~~~~~~~~~~~~~~~~~~~~~~~~~~~other~cases,~~~~~~~~~\\[2pt]
       \end{array}
    \right.
\end{align}
where $(\omega_{s,t})$ denotes the connection matrix of Levi-Civita connection $\nabla^L$.
\end{lem}
\begin{lem}{\rm \cite{Wa3}}
\begin{align}
\Gamma_{st}^k(x_0)&=\left\{
       \begin{array}{c}
        \Gamma^n_{ii}(x_0)=\frac{1}{2}h'(0),~~~~~~~~~~ ~~~~~~~~~~~{\rm if }~s=t=i,k=n,i<n; \\[2pt]
        \Gamma^i_{ni}(x_0)=-\frac{1}{2}h'(0),~~~~~~~~~~~~~~~~~~~{\rm if }~s=n,t=i,k=i,i<n;\\[2pt]
        \Gamma^i_{in}(x_0)=-\frac{1}{2}h'(0),~~~~~~~~~~~~~~~~~~~{\rm if }~s=i,t=n,k=i,i<n,\\[2pt]
        \Gamma_{st}^i(x_0)=0,~~~~~~~~~~~~~~~~~~~~~~~~~~~other~cases.~~~~~~~~~
       \end{array}
    \right.
\end{align}
\end{lem}
\indent Similar to (3.9) and (3.10), we firstly compute
\begin{equation}
\widetilde{{\rm Wres}}[\pi^+{{D}_{A}}^{-1}\circ\pi^+({D}^*_{A})^{-1}]=\int_M\int_{|\xi|=1}{\rm
trace}_{\wedge^*T^*M}[\sigma_{-4}(({{D}^*_{A}}{{D}_{A}})^{-1})]\sigma(\xi)dx+\int_{\partial M}\Psi,
\end{equation}
where
\begin{align}
\Psi &=\int_{|\xi'|=1}\int^{+\infty}_{-\infty}\sum^{\infty}_{j, k=0}\sum\frac{(-i)^{|\alpha|+j+k+1}}{\alpha!(j+k+1)!}
\times {\rm trace}_{\wedge^*T^*M}[\partial^j_{x_n}\partial^\alpha_{\xi'}\partial^k_{\xi_n}\sigma^+_{r}({{D}_{A}})^{-1}\\
&(x',0,\xi',\xi_n)\times\partial^\alpha_{x'}\partial^{j+1}_{\xi_n}\partial^k_{x_n}\sigma_{l}(({D}^*_{A})^{-1})(x',0,\xi',\xi_n)]d\xi_n\sigma(\xi')dx',\nonumber
\end{align}
the sum is taken over $r+l-k-j-|\alpha|-1=-4, r\leq -1, l\leq-1$.\\

\indent Then we can compute the interior of $\widetilde{{\rm Wres}}[\pi^+{{D}_{A}}^{-1}\circ\pi^+({D}^*_{A})^{-1}]$,
\begin{align}
&\int_M\int_{|\xi|=1}{\rm trace}_{\wedge^*T^*M}[\sigma_{-4}(({{D}^*_{A}}{{D}_{A}})^{-1})]\sigma(\xi)dx=32\pi^2\int_{M}{\rm tr}\bigg(-\frac{1}{12}s+A^*A\\
&-\frac{1}{4}\sum_{j=1}^{n}Ac(\widetilde{e_j})Ac(\widetilde{e_j})-\frac{1}{4}\sum_{j=1}^{n}A^*c(\widetilde{e_j})A^*c(\widetilde{e_j})+\frac{1}{2}\sum_{j=1}^{n}\nabla^{\bigwedge^*T^*M}_{\widetilde{e_j}}(A^*)c(\widetilde{e_j})\nonumber\\
&-\frac{1}{2}\sum_{j=1}^{n}c(\widetilde{e_j})\nabla^{\bigwedge^*T^*M}_{\widetilde{e_j}}(A)\bigg)d{\rm Vol_{M}}.\nonumber
\end{align}

\indent Now we  need to compute $\int_{\partial M} \Psi$. Since, some operators have the following symbols.
\begin{lem} The following identities hold:
\begin{align}
\sigma_1({D}_{A})&=\sigma_1({D}^*_{A})=ic(\xi);\\
\sigma_0({D}_{A})&=
\frac{1}{4}\sum_{i,s,t}\omega_{s,t}(\widetilde{e_i})c(\widetilde{e_i})\overline{c}(\widetilde{e_s})\overline{c}(\widetilde{e_t})-\frac{1}{4}\sum_{i,s,t}\omega_{s,t}(\widetilde{e_i})c(\widetilde{e_i})c(\widetilde{e_s})c(\widetilde{e_t})+A; \nonumber\\
\sigma_0({D}^*_{A})&=
\frac{1}{4}\sum_{i,s,t}\omega_{s,t}(\widetilde{e_i})c(\widetilde{e_i})\overline{c}(\widetilde{e_s})\overline{c}(\widetilde{e_t})-\frac{1}{4}\sum_{i,s,t}\omega_{s,t}(\widetilde{e_i})c(\widetilde{e_i})c(\widetilde{e_s})c(\widetilde{e_t})+A^*.\nonumber
\end{align}
\end{lem}
\indent Write
 \begin{eqnarray}
D_x^{\alpha}&=(-i)^{|\alpha|}\partial_x^{\alpha};
~\sigma({D}_{A})=p_1+p_0;
~\sigma({{D}_{A}}^{-1})=\sum^{\infty}_{j=1}q_{-j}.
\end{eqnarray}

\indent By the composition formula of pseudodifferential operators, we have
\begin{align}
1=\sigma({D}_{A}\circ {{D}_{A}}^{-1})&=\sum_{\alpha}\frac{1}{\alpha!}\partial^{\alpha}_{\xi}[\sigma({D}_{A})]
{D}_x^{\alpha}[\sigma({{D}_{A}}^{-1})]\\
&=(p_1+p_0)(q_{-1}+q_{-2}+q_{-3}+\cdots)\nonumber\\
&~~~+\sum_j(\partial_{\xi_j}p_1+\partial_{\xi_j}p_0)(
D_{x_j}q_{-1}+D_{x_j}q_{-2}+D_{x_j}q_{-3}+\cdots)\nonumber\\
&=p_1q_{-1}+(p_1q_{-2}+p_0q_{-1}+\sum_j\partial_{\xi_j}p_1D_{x_j}q_{-1})+\cdots,\nonumber
\end{align}
so
\begin{equation}
q_{-1}=p_1^{-1};~q_{-2}=-p_1^{-1}[p_0p_1^{-1}+\sum_j\partial_{\xi_j}p_1D_{x_j}(p_1^{-1})].
\end{equation}
\begin{lem} The following identities hold:
\begin{align}
\sigma_{-1}({{D}_{A}}^{-1})&=\sigma_{-1}(({D}^*_{A})^{-1})=\frac{ic(\xi)}{|\xi|^2};\\
\sigma_{-2}({{D}_{A}}^{-1})&=\frac{c(\xi)\sigma_{0}({D}_{A})c(\xi)}{|\xi|^4}+\frac{c(\xi)}{|\xi|^6}\sum_jc(dx_j)
\Big[\partial_{x_j}(c(\xi))|\xi|^2-c(\xi)\partial_{x_j}(|\xi|^2)\Big] ;\nonumber\\
\sigma_{-2}(({D}^*_{A})^{-1})&=\frac{c(\xi)\sigma_{0}({D}^*_{A})c(\xi)}{|\xi|^4}+\frac{c(\xi)}{|\xi|^6}\sum_jc(dx_j)
\Big[\partial_{x_j}(c(\xi))|\xi|^2-c(\xi)\partial_{x_j}(|\xi|^2)\Big].\nonumber
\end{align}
\end{lem}
\indent We denote ${\rm tr}$ as shorthand of ${\rm trace}.$ When $n=4$, then ${\rm tr}_{\wedge^*T^*M}[{\rm \texttt{id}}]={\rm dim}(\wedge^*(4))=16,$ since the sum is taken over $
r+l-k-j-|\alpha|-1=-4,~~r\leq -1,l\leq-1,$ then we have the following five cases:\\

\noindent  {\bf case a)~I)}~$r=-1,~l=-1,~k=j=0,~|\alpha|=1$\\

\noindent By applying the formula shown in (3.18), we can calculate
\begin{equation}
\Psi_1=-\int_{|\xi'|=1}\int^{+\infty}_{-\infty}\sum_{|\alpha|=1}
 {\rm trace}[\partial^\alpha_{\xi'}\pi^+_{\xi_n}\sigma_{-1}({{D}_{A}}^{-1})\times
 \partial^\alpha_{x'}\partial_{\xi_n}\sigma_{-1}(({D}_{A}^*)^{-1})](x_0)d\xi_n\sigma(\xi')dx'.
\end{equation}

\noindent  {\bf case a)~II)}~$r=-1,~l=-1,~k=|\alpha|=0,~j=1$\\

\noindent By (3.18), we get
\begin{equation}
\Psi_2=-\frac{1}{2}\int_{|\xi'|=1}\int^{+\infty}_{-\infty} {\rm
trace} [\partial_{x_n}\pi^+_{\xi_n}\sigma_{-1}({{D}_{A}}^{-1})\times
\partial_{\xi_n}^2\sigma_{-1}(({D}_{A}^*)^{-1})](x_0)d\xi_n\sigma(\xi')dx'.
\end{equation}

\noindent  {\bf case a)~III)}~$r=-1,~l=-1,~j=|\alpha|=0,~k=1$\\

\noindent By (3.18), we calculate that
\begin{equation}
\Psi_3=-\frac{1}{2}\int_{|\xi'|=1}\int^{+\infty}_{-\infty}
{\rm trace} [\partial_{\xi_n}\pi^+_{\xi_n}\sigma_{-1}({{D}_{A}}^{-1})\times
\partial_{\xi_n}\partial_{x_n}\sigma_{-1}(({D}_{A}^*)^{-1})](x_0)d\xi_n\sigma(\xi')dx'.\\
\end{equation}
Similar to the formulae (2.17)-(2.31) in \cite{Wa3}, we have
\begin{equation}
\Psi_1+\Psi_2+\Psi_3=0.
\end{equation}
\\
\noindent  {\bf case b)}~$r=-2,~l=-1,~k=j=|\alpha|=0$\\

\noindent Similarly, we get
\begin{align}
\Psi_4&=-i\int_{|\xi'|=1}\int^{+\infty}_{-\infty}{\rm trace} [\pi^+_{\xi_n}\sigma_{-2}({{D}_{A}}^{-1})\times
\partial_{\xi_n}\sigma_{-1}(({D}_{A}^*)^{-1})](x_0)d\xi_n\sigma(\xi')dx'.
\end{align}
We first compute\\
\begin{align}
\sigma_{-2}({{D}_{A}}^{-1})(x_0)=\frac{c(\xi)\sigma_{0}({{D}_{A}})(x_0)c(\xi)}{|\xi|^4}+\frac{c(\xi)}{|\xi|^6}c(dx_n)
\left[\partial_{x_n}[c(\xi')](x_0)|\xi|^2-c(\xi)h'(0)|\xi|^2_{\partial
M}\right],
\end{align}
where
\begin{align}
\sigma_{0}({{D}_{A}})(x_0)=
\frac{1}{4}\sum_{i,s,t}\omega_{s,t}(\widetilde{e_i})(x_{0})c(\widetilde{e_i})\overline{c}(\widetilde{e_s})\overline{c}(\widetilde{e_t})-\frac{1}{4}\sum_{i,s,t}\omega_{s,t}(\widetilde{e_i})(x_{0})c(\widetilde{e_i})c(\widetilde{e_s})c(\widetilde{e_t})+A.
\end{align}
We denote
\begin{align}
b_0^1(x_0)&=\frac{1}{4}\sum_{i,s,t}\omega_{s,t}(\widetilde{e_i})(x_{0})c(\widetilde{e_i})\overline{c}(\widetilde{e_s})\overline{c}(\widetilde{e_t});\\
b_0^2(x_0)&=-\frac{1}{4}\sum_{i,s,t}\omega_{s,t}(\widetilde{e_i})(x_{0})c(\widetilde{e_i})c(\widetilde{e_s})c(\widetilde{e_t}).
\end{align}
Means that
\begin{align}
&\pi^+_{\xi_n}\sigma_{-2}({{D}_{A}}^{-1})(x_0)|_{|\xi'|=1}=\pi^+_{\xi_n}\Big(\frac{c(\xi)b_0^1(x_0)c(\xi)}{(1+\xi_n^2)^2}\Big)+\pi^+_{\xi_n}\Big(\frac{c(\xi)Ac(\xi)}{(1+\xi_n^2)^2}\Big)
\nonumber\\
&+\pi^+_{\xi_n}\Big(\frac{c(\xi)b_0^2(x_0)c(\xi)+c(\xi)c(dx_n)\partial_{x_n}[c(\xi')](x_0)}{(1+\xi_n^2)^2}-h'(0)\frac{c(\xi)c(dx_n)c(\xi)}{(1+\xi_n^{2})^3}\Big).\nonumber
\end{align}
Since
\begin{align}
b_0^1(x_0)c(dx_n)=-\frac{1}{4}h'(0)\sum_{i=1}^{n-1}c(\widetilde{e_i})\overline{c}(\widetilde{e_i})c(\widetilde{e_n})\overline{c}(\widetilde{e_n}),
\end{align}
then by the relation of the Clifford action and ${\rm tr}{AB}={\rm tr}{BA}$,  we have the equalities:\\
\begin{align}
{\rm tr}[c(\widetilde{e_i})\overline{c}(\widetilde{e_i})c(\widetilde{e_n})\overline{c}(\widetilde{e_n})]=0(i<n);~
{\rm tr}[b_0^1(x_0)c(dx_n)]=0.\nonumber
\end{align}
Therefore
\begin{align}
\partial_{\xi_n}\sigma_{-1}(({D}_{A}^*)^{-1})(x_0)|_{|\xi'|=1}=\partial_{\xi_n}q_{-1}(x_0)|_{|\xi'|=1}=i\left(\frac{c(dx_n)}{1+\xi_n^2}-\frac{2\xi_nc(\xi')+2\xi_n^2c(dx_n)}{(1+\xi_n^2)^2}\right).
\end{align}
Hence, we have
\begin{align}
&{\rm tr}\Big[\pi^+_{\xi_n}\bigg(\frac{c(\xi)b_0^1(x_0)c(\xi)}{(1+\xi_n^2)^2}\bigg)\times
\partial_{\xi_n}\sigma_{-1}(({D}_{A}^*)^{-1})(x_0)\Big]|_{|\xi'|=1}\\
&=\frac{1}{2(1+\xi_n^2)^2}{\rm tr}[b_0^1(x_0)c(\xi')]+\frac{i}{2(1+\xi_n^2)^2}{\rm tr}[b_0^1(x_0)c(dx_n)]\nonumber\\
&=\frac{1}{2(1+\xi_n^2)^2}{\rm tr}[b_0^1(x_0)c(\xi')].\nonumber
\end{align}
We note that $i<n,$  $\int_{|\xi'|=1}{\{\xi_{i_1}\cdot\cdot\cdot\xi_{i_{2d+1}}}\}\sigma(\xi')=0,$ so ${\rm tr}[c(\xi')c(dx_n)]$ has no contribution for computing case b),
\begin{align}
&-i\int_{|\xi'|=1}\int^{+\infty}_{-\infty}{\rm trace}\Big[\pi^+_{\xi_n}\Big(\frac{c(\xi)b_0^1(x_0)c(\xi)}{(1+\xi_n^2)^2}\Big)\times
\partial_{\xi_n}\sigma_{-1}(({D}_{A}^*)^{-1})\Big](x_0)d\xi_n\sigma(\xi')dx'=0.
\end{align}
Since
\begin{align}
&{\rm tr}\Big[\pi^+_{\xi_n}\Big(\frac{c(\xi)Ac(\xi)}{(1+\xi_n^2)^2}\Big)\times
\partial_{\xi_n}\sigma_{-1}(({D}_{A}^*)^{-1})(x_0)\Big]|_{|\xi'|=1}\\
&=\frac{1}{2(1+\xi_n^2)^2}{\rm tr}[Ac(\xi')]+\frac{i}{2(1+\xi_n^2)^2}{\rm tr}[Ac(dx_n)],\nonumber
\end{align}
then, we have
\begin{align}
&-i\int_{|\xi'|=1}\int^{+\infty}_{-\infty}{\rm trace}\Big[\pi^+_{\xi_n}\Big(\frac{c(\xi)Ac(\xi)}{(1+\xi_n^2)^2}\Big)\times
\partial_{\xi_n}\sigma_{-1}(({D}_{A}^*)^{-1})\Big](x_0)d\xi_n\sigma(\xi')dx'\\
&=-i\int_{|\xi'|=1}\int^{+\infty}_{-\infty}\frac{1}{2(1+\xi_n^2)^2}{\rm tr}[Ac(\xi')]d\xi_n\sigma(\xi')dx'\nonumber\\
&-i\int_{|\xi'|=1}\int^{+\infty}_{-\infty}\frac{i}{2(1+\xi_n^2)^2}{\rm tr}[Ac(dx_n)]d\xi_n\sigma(\xi')dx'\nonumber\\
&=\frac{\Omega_3}{2}{\rm tr}[Ac(dx_n)]\int_{\Gamma^+}\frac{1}{(\xi_n-i)^2(\xi_n+i)^2}dx'\nonumber\\
&=\frac{\Omega_3}{2}{\rm tr}[Ac(dx_n)]2\pi i[\frac{1}{(\xi_n+i)^2}]^{(1)}|_{\xi_n=i}dx'\nonumber\\
&=\frac{\pi}{4}\Omega_3{\rm tr}[Ac(dx_n)]dx'.\nonumber
\end{align}
We have
\begin{equation}
\pi^+_{\xi_n}\Big(\frac{c(\xi)b_0^2(x_0)c(\xi)+c(\xi)c(dx_n)\partial_{x_n}[c(\xi')](x_0)}{(1+\xi_n^2)^2}\Big)-h'(0)\pi^+_{\xi_n}\Big(\frac{c(\xi)c(dx_n)c(\xi)}{(1+\xi_n)^3}\Big):= B_1-B_2,
\end{equation}
where
\begin{align}
B_1&=\frac{-1}{4(\xi_n-i)^2}[(2+i\xi_n)c(\xi')b_0^2(x_0)c(\xi')+i\xi_nc(dx_n)b_0^2(x_0)c(dx_n)-i\partial_{x_n}c(\xi')\\
&+(2+i\xi_n)c(\xi')c(dx_n)\partial_{x_n}c(\xi')+ic(dx_n)b_0^2(x_0)c(\xi')+ic(\xi')b_0^2(x_0)c(dx_n)]\nonumber
\end{align}
and
\begin{align}
B_2&=\frac{h'(0)}{2}\left(\frac{c(dx_n)}{4i(\xi_n-i)}+\frac{c(dx_n)-ic(\xi')}{8(\xi_n-i)^2}
+\frac{3\xi_n-7i}{8(\xi_n-i)^3}(ic(\xi')-c(dx_n))\right).
\end{align}
A simple calculation shows that
\begin{align}
{\rm tr }[B_2\times\partial_{\xi_n}\sigma_{-1}(({D}_{A}^*)^{-1})(x_0)]|_{|\xi'|=1}
&=\frac{i}{2}h'(0)\frac{-i\xi_n^2-\xi_n+4i}{4(\xi_n-i)^3(\xi_n+i)^2}{\rm tr}[\texttt{id}]\\
&=8ih'(0)\frac{-i\xi_n^2-\xi_n+4i}{4(\xi_n-i)^3(\xi_n+i)^2}.\nonumber
\end{align}
Similarly, we have
\begin{align}
{\rm tr }[B_1\times\partial_{\xi_n}\sigma_{-1}(({D}_{A}^*)^{-1})(x_0)]|_{|\xi'|=1}=
\frac{-8ic_0}{(1+\xi_n^2)^2}+2h'(0)\frac{\xi_n^2-i\xi_n-2}{(\xi_n-i)(1+\xi_n^2)^2},
\end{align}
where $b_0^2=c_0c(dx_n)$ and $c_0=-\frac{3}{4}h'(0)$.\\
By (3.44) and (3.45), we have
\begin{align}
&-i\int_{|\xi'|=1}\int^{+\infty}_{-\infty}{\rm trace} [(B_1-B_2)\times
\partial_{\xi_n}\sigma_{-1}(({{D}_{A}^*})^{-1})](x_0)d\xi_n\sigma(\xi')dx'\\
&=-\Omega_3\int_{\Gamma^+}\frac{8c_0(\xi_n-i)+4ih'(0)}{(\xi_n-i)^3(\xi_n+i)^2}d\xi_ndx'\nonumber\\
&=\frac{9}{2}\pi h'(0)\Omega_3dx'.\nonumber
\end{align}
Hence, we have
\begin{align}
\Psi_4=\frac{9}{2}\pi h'(0)\Omega_3dx'+\frac{\pi}{4}\Omega_3{\rm tr}[Ac(dx_n)]dx'.
\end{align}

\noindent {\bf  case c)}~$r=-1,~l=-2,~k=j=|\alpha|=0$\\

\noindent Using (3.18), we get
\begin{align}
\Psi_5=-i\int_{|\xi'|=1}\int^{+\infty}_{-\infty}{\rm trace} [\pi^+_{\xi_n}\sigma_{-1}({{D}_{A}}^{-1})\times
\partial_{\xi_n}\sigma_{-2}(({D}^*_{A})^{-1})](x_0)d\xi_n\sigma(\xi')dx'.
\end{align}
Considering (3.5) and (3.6), we have
\begin{align}
\pi^+_{\xi_n}\sigma_{-1}({{D}_{A}}^{-1})(x_0)|_{|\xi'|=1}=\frac{c(\xi')+ic(dx_n)}{2(\xi_n-i)}.
\end{align}
Since
\begin{equation}
\sigma_{-2}(({D}^*_{A})^{-1})(x_0)=\frac{c(\xi)\sigma_{0}({D}^*_{A})(x_0)c(\xi)}{|\xi|^4}+\frac{c(\xi)}{|\xi|^6}c(dx_n)
\left[\partial_{x_n}[c(\xi')](x_0)|\xi|^2-c(\xi)h'(0)|\xi|^2_{\partial_M}\right],
\end{equation}
where
\begin{align}
\sigma_{0}({D}^*_{A})(x_0)&=
\frac{1}{4}\sum_{i,s,t}\omega_{s,t}(\widetilde{e_i})(x_{0})c(\widetilde{e_i})\overline{c}(\widetilde{e_s})\overline{c}(\widetilde{e_t})-\frac{1}{4}\sum_{i,s,t}\omega_{s,t}(\widetilde{e_i})(x_{0})c(\widetilde{e_i})c(\widetilde{e_s})c(\widetilde{e_t})+A^*\\
&=b_0^1(x_0)+b_0^2(x_0)+A^*,\nonumber
\end{align}
then
\begin{align}
\partial_{\xi_n}\sigma_{-2}(({D}^*_{A})^{-1})(x_0)|_{|\xi'|=1}&=
\partial_{\xi_n}\bigg(\frac{c(\xi)(b_0^1(x_0)+b_0^2(x_0)+A^*)c(\xi)}{|\xi|^4}\\
&+\frac{c(\xi)}{|\xi|^6}c(dx_n)\left[\partial_{x_n}[c(\xi')](x_0)|\xi|^2-c(\xi)h'(0)|\xi|^2_{\partial_M}\right]\bigg)\nonumber\\
&=\partial_{\xi_n}\bigg(\frac{c(\xi)b_0^2(x_0)c(\xi)}{|\xi|^4}+\frac{c(\xi)}{|\xi|^6}c(dx_n)\left[\partial_{x_n}[c(\xi')](x_0)|\xi|^2-c(\xi)h'(0)\right]\bigg)\nonumber\\
&+\partial_{\xi_n}\bigg(\frac{c(\xi)b_0^1(x_0)c(\xi)}{|\xi|^4}\bigg)+\partial_{\xi_n}\bigg(\frac{c(\xi)A^*c(\xi)}{|\xi|^4}\bigg).\nonumber
\end{align}
By computation, we have
\begin{align}
\partial_{\xi_n}\bigg(\frac{c(\xi)b_0^1(x_0)c(\xi)}{|\xi|^4}\bigg)&=
\frac{c(dx_n)b_0^1(x_0)c(\xi)}{|\xi|^4}+\frac{c(\xi)b_0^1(x_0)c(dx_n)}{|\xi|^4}-\frac{4\xi_nc(\xi)b_0^1(x_0)c(\xi)}{|\xi|^6};\\
\partial_{\xi_n}\bigg(\frac{c(\xi)A^*c(\xi)}{|\xi|^4}\bigg)&=
\frac{c(dx_n)A^*c(\xi)}{|\xi|^4}+\frac{c(\xi)A^*c(dx_n)}{|\xi|^4}-\frac{4\xi_nc(\xi)A^*c(\xi)}{|\xi|^6}.
\end{align}
For the sake of convenience in writing, we denote
\begin{align}
q_{-2}^{1}=\frac{c(\xi)b_0^2(x_0)c(\xi)}{|\xi|^4}+\frac{c(\xi)}{|\xi|^6}c(dx_n)\left[\partial_{x_n}[c(\xi')](x_0)|\xi|^2-c(\xi)h'(0)\right],
\end{align}
then
\begin{align}
\partial_{\xi_n}(q_{-2}^{1})&=\frac{1}{(1+\xi_n^2)^3}[(2\xi_n-2\xi_n^3)c(dx_n)b_0^2(x_0)c(dx_n)+(1-3\xi_n^2)c(dx_n)b_0^2(x_0)c(\xi')\\
&+(1-3\xi_n^2)c(\xi')b_0^2(x_0)c(dx_n)-4\xi_nc(\xi')b_0^2(x_0)c(\xi')+(3\xi_n^2-1)\partial_{x_n}c(\xi')\nonumber\\
&-4\xi_nc(\xi')c(dx_n)\partial_{x_n}c(\xi')+2h'(0)c(\xi')+2h'(0)\xi_nc(dx_n)]\nonumber\\
&+6\xi_nh'(0)\frac{c(\xi)c(dx_n)c(\xi)}{(1+\xi^2_n)^4}.\nonumber
\end{align}
By (3.48) and (3.52), we have
\begin{align}
&{\rm tr}\bigg[\pi^+_{\xi_n}\sigma_{-1}({{D}_{A}}^{-1})(x_0)\times
\partial_{\xi_n}\bigg(\frac{c(\xi)b_0^1(x_0)c(\xi)}{|\xi|^4}\bigg)\bigg]|_{|\xi'|=1}\\
&=\frac{-1}{2(\xi_n-i)(\xi_n+i)^3}{\rm tr}[b_0^1(x_0)c(\xi')]+\frac{i}{2(\xi_n-i)(\xi_n+i)^3}{\rm tr}[b_0^1(x_0)c(dx_n)]\nonumber\\
&=\frac{-1}{2(\xi_n-i)(\xi_n+i)^3}{\rm tr}[b_0^1(x_0)c(\xi')],\nonumber
\end{align}
it is shown that
\begin{align}
&-i\int_{|\xi'|=1}\int^{+\infty}_{-\infty}{\rm trace}\bigg[\pi^+_{\xi_n}\sigma_{-1}({{D}_{A}}^{-1})\times
\partial_{\xi_n}\bigg(\frac{c(\xi)b_0^1c(\xi)}{|\xi|^4}\bigg)\bigg](x_0)d\xi_n\sigma(\xi')dx'\\
&=-i\int_{|\xi'|=1}\int^{+\infty}_{-\infty}\frac{-1}{2(\xi_n-i)(\xi_n+i)^3}{\rm tr}[b_0^1(x_0)c(\xi')]d\xi_n\sigma(\xi')dx'\nonumber\\
&=0.\nonumber
\end{align}
Similarly to (3.56), we have
\begin{align}
&{\rm tr}\bigg[\pi^+_{\xi_n}\sigma_{-1}({{D}_{A}}^{-1})(x_0)\times
\partial_{\xi_n}\bigg(\frac{c(\xi)A^*c(\xi)}{|\xi|^4}\bigg)\bigg]|_{|\xi'|=1}\\
&=\frac{-1}{2(\xi_n-i)(\xi_n+i)^3}{\rm tr}[A^*c(\xi')]+\frac{i}{2(\xi_n-i)(\xi_n+i)^3}{\rm tr}[A^*c(dx_n)],\nonumber
\end{align}
then
\begin{align}
&-i\int_{|\xi'|=1}\int^{+\infty}_{-\infty}{\rm trace}\bigg[\pi^+_{\xi_n}\sigma_{-1}({{D}_{A}}^{-1})\times
\partial_{\xi_n}\bigg(\frac{c(\xi)A^*c(\xi)}{|\xi|^4}\bigg)\bigg](x_0)d\xi_n\sigma(\xi')dx'\\
&=-i\int_{|\xi'|=1}\int^{+\infty}_{-\infty}\frac{-1}{2(\xi_n-i)(\xi_n+i)^3}{\rm tr}[A^*c(\xi')]d\xi_n\sigma(\xi')dx'\nonumber\\
&-i\int_{|\xi'|=1}\int^{+\infty}_{-\infty}\frac{i}{2(\xi_n-i)(\xi_n+i)^3}{\rm tr}[A^*c(dx_n)]d\xi_n\sigma(\xi')dx'\nonumber\\
&=\frac{\Omega_3}{2}{\rm tr}[A^*c(dx_n)]\int_{\Gamma^+}\frac{1}{(\xi_n-i)(\xi_n+i)^3}d\xi_ndx'\nonumber\\
&=\frac{\Omega_3}{2}{\rm tr}[A^*c(dx_n)]2\pi i[\frac{1}{(\xi_n+i)^3}]^{(1)}|_{\xi_n=i}dx'\nonumber\\
&=-\frac{\pi}{4}\Omega_3{\rm tr}[A^*c(dx_n)]dx'.\nonumber
\end{align}
Observing (3.48) and (3.55), we have
\begin{align}
{\rm tr}[\pi^+_{\xi_n}\sigma_{-1}({{D}_{A}}^{-1})\times
\partial_{\xi_n}(q^1_{-2})](x_0)|_{|\xi'|=1}
=\frac{12h'(0)(i\xi^2_n+\xi_n-2i)}{(\xi-i)^3(\xi+i)^3}
+\frac{48h'(0)i\xi_n}{(\xi-i)^3(\xi+i)^4}.
\end{align}
By $\int_{|\xi'|=1}{\{\xi_{i_1}\cdot\cdot\cdot\xi_{i_{2d+1}}}\}\sigma(\xi')=0$ and (3.60), we have
\begin{align}
&-i\int_{|\xi'|=1}\int^{+\infty}_{-\infty}{\rm tr}[\pi^+_{\xi_n}\sigma_{-1}({{D}_{A}}^{-1})\times
\partial_{\xi_n}(q^1_{-2})](x_0)d\xi_n\sigma(\xi')dx'\\
&=-i\int_{|\xi'|=1}\int^{+\infty}_{-\infty}\frac{12h'(0)(i\xi^2_n+\xi_n-2i)}{(\xi-i)^3(\xi+i)^3}+\frac{48h'(0)i\xi_n}{(\xi-i)^3(\xi+i)^4}d\xi_n\sigma(\xi')dx'\nonumber\\
&=-i\Omega_3\int_{\Gamma^+}\frac{12h'(0)(i\xi^2_n+\xi_n-2i)}{(\xi-i)^3(\xi+i)^3}+\frac{48h'(0)i\xi_n}{(\xi-i)^3(\xi+i)^4}d\xi_ndx'\nonumber\\
&=-i\Omega_3\frac{2\pi i}{2!}[\frac{12h'(0)(i\xi^2_n+\xi_n-2i)}{(\xi+i)^3}]^{(2)}|_{\xi_n=i}dx'-i\Omega_3\frac{2\pi i}{2!}[\frac{48h'(0)i\xi_n}{(\xi+i)^4}]^{(2)}|_{\xi_n=i}dx'\nonumber\\
&=-\frac{9}{2}\pi h'(0)\Omega_3dx'.\nonumber
\end{align}
Then,
\begin{align}
\Psi_5=-\frac{9}{2}\pi h'(0)\Omega_3dx'-\frac{\pi}{4}\Omega_3{\rm tr}[A^*c(dx_n)]dx'.
\end{align}
In summary,
\begin{align}
\Psi=\Psi_1+\Psi_2+\Psi_3+\Psi_4+\Psi_5=\frac{\pi}{4}\Omega_3{\rm tr}[Ac(dx_n)]dx'-\frac{\pi}{4}\Omega_3{\rm tr}[A^*c(dx_n)]dx'.
\end{align}

Applying (3.17), (3.19) and (3.63), we can assert that:
\begin{thm}
Let $M$ be a $4$-dimensional oriented
compact manifolds with the boundary $\partial M$ and the metric
$g^M$ as above, ${{D}_{A}}$ be the perturbation of the de Rham Hodge operator on $\widetilde{M}$, then
\begin{align}
&\widetilde{{\rm Wres}}[\pi^+{{D}_{A}}^{-1}\circ\pi^+({D}^*_{A})^{-1}]=32\pi^2\int_{M}{\rm tr}\bigg(-\frac{1}{12}s+A^*A-\frac{1}{4}\sum_{j=1}^{n}Ac(\widetilde{e_j})Ac(\widetilde{e_j})\\
&-\frac{1}{4}\sum_{j=1}^{n}A^*c(\widetilde{e_j})A^*c(\widetilde{e_j})+\frac{1}{2}\sum_{j=1}^{n}\nabla^{\bigwedge^*T^*M}_{\widetilde{e_j}}(A^*)c(\widetilde{e_j})-\frac{1}{2}\sum_{j=1}^{n}c(\widetilde{e_j})\nabla^{\bigwedge^*T^*M}_{\widetilde{e_j}}(A)\bigg)d{\rm Vol_{M}}\nonumber\\
&+\int_{\partial M}\frac{\pi}{4}\Omega_3{\rm tr}[Ac(dx_n)]d{\rm Vol_{\partial M}}-\int_{\partial M}\frac{\pi}{4}\Omega_3{\rm tr}[A^*c(dx_n)]d{\rm Vol_{\partial M}},\nonumber
\end{align}
where $s$ is the scalar curvature.
\end{thm}

When $A=c(X),$ then we have
\begin{align}
&\int_M\int_{|\xi|=1}{\rm
trace}_{\wedge^*T^*M}[\sigma_{-4}(({{D}^*_{A}}{{D}_{A}})^{-1})]\sigma(\xi)dx\\
&=512\pi^2\int_{M}\bigg(-\frac{1}{12}s+2|X|^2+\sum_{j=1}^{n}g(\nabla^{TM}_{\widetilde{e_j}}X, \widetilde{e_j})\bigg)d{\rm Vol_{M}}\nonumber
\end{align}
and
\begin{align}
\int_{\partial M} \Psi=-\int_{\partial M}8\pi\Omega_3 g(\partial{x_n}, X)d{\rm Vol_{\partial M}}.
\end{align}
We can immediately state the following corollary:
\begin{cor}
Let $M$ be a $4$-dimensional oriented
compact manifolds with the boundary $\partial M$ and the metric
$g^M$ as above, and let $A=c(X),$ then
\begin{align}
&\widetilde{{\rm Wres}}[\pi^+{{D}_{A}}^{-1}\circ\pi^+({D}^*_{A})^{-1}]=512\pi^2\int_{M}\bigg(-\frac{1}{12}s+2|X|^2+\sum_{j=1}^{n}g(\nabla^{TM}_{\widetilde{e_j}}X, \widetilde{e_j})\bigg)d{\rm Vol_{M}}\\
&-\int_{\partial M}8\pi\Omega_3 g(\partial{x_n}, X)d{\rm Vol_{\partial M}},\nonumber
\end{align}
where $s$ is the scalar curvature.
\end{cor}

When $A=\overline{c}(X),$ we can get
\begin{align}
&\int_M\int_{|\xi|=1}{\rm
trace}_{\wedge^*T^*M}[\sigma_{-4}(({{D}^*_{A}}{{D}_{A}})^{-1})]\sigma(\xi)dx=512\pi^2\int_{M}\bigg(-\frac{1}{12}s-|X|^2\bigg)d{\rm Vol_{M}},
\end{align}
and
\begin{align}
\int_{\partial M} \Psi=0.
\end{align}
Now, we compute $\widetilde{{\rm Wres}}[\pi^+{{D}_{A}}^{-1}\circ\pi^+({D}^*_{A})^{-1}].$
\begin{cor}
Let $M$ be a $4$-dimensional oriented
compact manifolds with the boundary $\partial M$ and the metric
$g^M$ as above, and let $A=\overline{c}(X),$ then
\begin{align}
\widetilde{{\rm Wres}}[\pi^+{{D}_{A}}^{-1}\circ\pi^+({D}^*_{A})^{-1}]&=512\pi^2\int_{M}\bigg(-\frac{1}{12}s-|X|^2\bigg)d{\rm Vol_{M}},
\end{align}
where $s$ is the scalar curvature.
\end{cor}

When $A=c(X)c(Y),$ then  $A^*=c(Y)c(X).$\\
By computation, we have
\begin{align}
&{\rm tr}[A^*A]={\rm tr}[(-|X|^2)(-|Y|^2)]=|X|^2|Y|^2{\rm tr}[\texttt{id}],\\
&{\rm tr}[c(dx_n)A]=0, {\rm tr}[c(dx_n)A^*]=0,\nonumber
\end{align}
\begin{align}
&\sum_{j=1}^{n}{\rm tr}[Ac(\widetilde{e_j})Ac(\widetilde{e_j})]=\sum_{j=1}^{n}{\rm tr}[c(X)c(Y)c(\widetilde{e_j})c(X)c(Y)c(\widetilde{e_j})]\\
&=-\sum_{j=1}^{n}{\rm tr}[c(X)c(Y)c(X)c(\widetilde{e_j})c(Y)c(\widetilde{e_j})]-2\sum_{j=1}^{n}g(\widetilde{e_j},X){\rm tr}[c(X)c(Y)c(Y)c(\widetilde{e_j})]\nonumber\\
&=-\sum_{j=1}^{n}{\rm tr}[c(X)c(Y)c(X)c(\widetilde{e_j})c(Y)c(\widetilde{e_j})]-2|X|^2|Y|^2{\rm tr}[\texttt{id}]\nonumber\\
&=(n-4)|X|^2|Y|^2{\rm tr}[\texttt{id}]+(4-2n)g(X,Y)^2{\rm tr}[\texttt{id}]\nonumber\\
&\sum_{j=1}^{n}{\rm tr}[A^*c(\widetilde{e_j})A^*c(\widetilde{e_j})]=\sum_{j=1}^{n}{\rm tr}[c(Y)c(X)c(\widetilde{e_j})c(Y)c(X)c(\widetilde{e_j})]\\
&=-\sum_{j=1}^{n}{\rm tr}[c(Y)c(X)c(Y)c(\widetilde{e_j})c(X)c(\widetilde{e_j})]-2\sum_{j=1}^{n}g(\widetilde{e_j},Y){\rm tr}[c(Y)c(X)c(X)c(\widetilde{e_j})]\nonumber\\
&=-\sum_{j=1}^{n}{\rm tr}[c(Y)c(X)c(Y)c(\widetilde{e_j})c(X)c(\widetilde{e_j})]-2|X|^2|Y|^2{\rm tr}[\texttt{id}]\nonumber\\
&=(n-4)|X|^2|Y|^2{\rm tr}[\texttt{id}]+(4-2n)g(X,Y)^2{\rm tr}[\texttt{id}]\nonumber\\
&\sum_{j=1}^{n}{\rm tr}[\nabla^{\bigwedge^*T^*M}_{\widetilde{e_j}}(A^*)c(\widetilde{e_j})]=\sum_{j=1}^{n}{\rm tr}[\nabla^{\bigwedge^*T^*M}_{\widetilde{e_j}}\big(c(Y)c(X)\big)c(\widetilde{e_j})]\\
&=\sum_{j=1}^{n}{\rm tr}[\nabla^{\bigwedge^*T^*M}_{\widetilde{e_j}}\big(c(Y)\big)c(X)c(\widetilde{e_j})+c(Y)\nabla^{\bigwedge^*T^*M}_{\widetilde{e_j}}\big(c(X)\big)c(\widetilde{e_j})]\nonumber\\
&=\sum_{j=1}^{n}{\rm tr}[c(\nabla^{TM}_{\widetilde{e_j}}Y)c(X)c(\widetilde{e_j})+c(Y)c(\nabla^{TM}_{\widetilde{e_j}}X)c(\widetilde{e_j})]=0,\nonumber\\
&\sum_{j=1}^{n}{\rm tr}[c(\widetilde{e_j})\nabla^{\bigwedge^*T^*M}_{\widetilde{e_j}}(A)]=\sum_{j=1}^{n}{\rm tr}[c(\widetilde{e_j})\nabla^{\bigwedge^*T^*M}_{\widetilde{e_j}}\big(c(X)c(Y)\big)]\\
&=\sum_{j=1}^{n}{\rm tr}[c(\widetilde{e_j})\nabla^{\bigwedge^*T^*M}_{\widetilde{e_j}}\big(c(X)\big)c(Y)+c(\widetilde{e_j})c(X)\nabla^{\bigwedge^*T^*M}_{\widetilde{e_j}}\big(c(Y)\big)]\nonumber\\
&=\sum_{j=1}^{n}{\rm tr}[c(\widetilde{e_j})c(\nabla^{TM}_{\widetilde{e_j}}X)c(Y)+c(\widetilde{e_j})c(X)c(\nabla^{TM}_{\widetilde{e_j}}Y)]=0.\nonumber
\end{align}
By applying the formulae shown in (3.19), we can calculate
\begin{align}
&\int_M\int_{|\xi|=1}{\rm
trace}_{\wedge^*T^*M}[\sigma_{-4}(({{D}^*_{A}}{{D}_{A}})^{-1})]\sigma(\xi)dx\\
&=32\pi^2\int_{M}{\rm tr}\bigg(-\frac{1}{12}s+A^*A-\frac{1}{4}\sum_{j=1}^{n}Ac(\widetilde{e_j})Ac(\widetilde{e_j})-\frac{1}{4}\sum_{j=1}^{n}A^*c(\widetilde{e_j})A^*c(\widetilde{e_j})\nonumber\\
&+\frac{1}{2}\sum_{j=1}^{n}\nabla^{\bigwedge^*T^*M}_{\widetilde{e_j}}(A^*)c(\widetilde{e_j})-\frac{1}{2}\sum_{j=1}^{n}c(\widetilde{e_j})\nabla^{\bigwedge^*T^*M}_{\widetilde{e_j}}(A)\bigg)d{\rm Vol_{M}}\nonumber\\
&=32\pi^2\int_{M}\bigg(-\frac{1}{12}s+|X|^2|Y|^2-\frac{1}{4}(-4)g(X,Y)^2-\frac{1}{4}(-4)g(X,Y)^2\bigg){\rm tr}[\texttt{id}]d{\rm Vol_{M}}\nonumber\\
&=512\pi^2\int_{M}\bigg(-\frac{1}{12}s+|X|^2|Y|^2+2g(X,Y)^2\bigg)d{\rm Vol_{M}},\nonumber
\end{align}
and
\begin{align}
\int_{\partial M} \Psi=0.
\end{align}
We can claim the following corollary:
\begin{cor}
Let $M$ be a $4$-dimensional oriented
compact manifolds with the boundary $\partial M$ and the metric
$g^M$ as above, and let $A=c(X)c(Y),$ then
\begin{align}
&\widetilde{{\rm Wres}}[\pi^+{{D}_{A}}^{-1}\circ\pi^+({D}^*_{A})^{-1}]=512\pi^2\int_{M}\bigg(-\frac{1}{12}s+|X|^2|Y|^2+2g(X,Y)^2\bigg)d{\rm Vol_{M}},
\end{align}
where $s$ is the scalar curvature.
\end{cor}

When $A=c(X)\overline{c}(Y),$ similar to (3.78), we can get:
\begin{cor}
Let $M$ be a $4$-dimensional oriented
compact manifolds with the boundary $\partial M$ and the metric
$g^M$ as above, and let $A=c(X)\overline{c}(Y),$ then
\begin{align}
&\widetilde{{\rm Wres}}[\pi^+{{D}_{A}}^{-1}\circ\pi^+({D}^*_{A})^{-1}]=512\pi^2\int_{M}\bigg(-\frac{1}{12}s+2|X|^2|Y|^2\bigg)d{\rm Vol_{M}},
\end{align}
where $s$ is the scalar curvature.
\end{cor}

When $A=\overline{c}(X)\overline{c}(Y),$ similar to (3.78), we can get the following corollary:
\begin{cor}
Let $M$ be a $4$-dimensional oriented
compact manifolds with the boundary $\partial M$ and the metric
$g^M$ as above, and let $A=\overline{c}(X)\overline{c}(Y),$ then
\begin{align}
&\widetilde{{\rm Wres}}[\pi^+{{D}_{A}}^{-1}\circ\pi^+({D}^*_{A})^{-1}]=512\pi^2\int_{M}\bigg(-\frac{1}{12}s-|X|^2|Y|^2+4g(X,Y)^2\bigg)d{\rm Vol_{M}},
\end{align}
where $s$ is the scalar curvature.
\end{cor}

When $A=c(X)c(Y)c(Z),$ then  $A^*=-c(Z)c(Y)c(X).$\\
By computation, we have
\begin{align}
&{\rm tr}[A^*A]=-{\rm tr}[(-|X|^2)(-|Y|^2)(-|Z|^2)]=|X|^2|Y|^2|Z|^2{\rm tr}[\texttt{id}],\\
&\sum_{j=1}^{n}{\rm tr}[Ac(\widetilde{e_j})Ac(\widetilde{e_j})]=\sum_{j=1}^{n}{\rm tr}[c(X)c(Y)c(Z)c(\widetilde{e_j})c(X)c(Y)c(Z)c(\widetilde{e_j})]\\
&=(n-6)|X|^2|Y|^2|Z|^2{\rm tr}[\texttt{id}]+(8-2n)|X|^2g(Y,Z)^2{\rm tr}[\texttt{id}]+(8-2n)|Y|^2g(X,Z)^2{\rm tr}[\texttt{id}]\nonumber\\
&+(8-2n)|Z|^2g(X,Y)^2{\rm tr}[\texttt{id}]+(4n-16)g(X,Y)g(X,Z)g(Y,Z){\rm tr}[\texttt{id}]\nonumber\\
&\sum_{j=1}^{n}{\rm tr}[A^*c(\widetilde{e_j})A^*c(\widetilde{e_j})]=\sum_{j=1}^{n}{\rm tr}[(-c(Z)c(Y)c(X))c(\widetilde{e_j})(-c(Z)c(Y)c(X))c(\widetilde{e_j})]\\
&=(n-6)|X|^2|Y|^2|Z|^2{\rm tr}[\texttt{id}]+(8-2n)|X|^2g(Y,Z)^2{\rm tr}[\texttt{id}]+(8-2n)|Y|^2g(X,Z)^2{\rm tr}[\texttt{id}]\nonumber\\
&+(8-2n)|Z|^2g(X,Y)^2{\rm tr}[\texttt{id}]+(4n-16)g(X,Y)g(X,Z)g(Y,Z){\rm tr}[\texttt{id}]\nonumber\\
&\sum_{j=1}^{n}{\rm tr}[\nabla^{\bigwedge^*T^*M}_{\widetilde{e_j}}(A^*)c(\widetilde{e_j})]=\sum_{j=1}^{n}{\rm tr}[\nabla^{\bigwedge^*T^*M}_{\widetilde{e_j}}\big(-c(Z)c(Y)c(X)\big)c(\widetilde{e_j})]\\
&=\big[-\sum_{j=1}^{n}g(\nabla^{TM}_{\widetilde{e_j}}Z, \widetilde{e_j})g(X, Y)+\sum_{j=1}^{n}g(Y, \widetilde{e_j})g(X, \nabla^{TM}_{\widetilde{e_j}}Z)-\sum_{j=1}^{n}g(X, \widetilde{e_j})g(Y, \nabla^{TM}_{\widetilde{e_j}}Z)\nonumber\\
&-\sum_{j=1}^{n}g(Z, \widetilde{e_j})g(X,\nabla^{TM}_{\widetilde{e_j}}Y)+\sum_{j=1}^{n}g(\nabla^{TM}_{\widetilde{e_j}}Y, \widetilde{e_j})g(X, Z)-\sum_{j=1}^{n}g(X, \widetilde{e_j})g(\nabla^{TM}_{\widetilde{e_j}}Y, Z)\nonumber\\
&-\sum_{j=1}^{n}g(Z, \widetilde{e_j})g(\nabla^{TM}_{\widetilde{e_j}}X, Y)+\sum_{j=1}^{n}g(Y, \widetilde{e_j})g(\nabla^{TM}_{\widetilde{e_j}}X, Z)-\sum_{j=1}^{n}g(\nabla^{TM}_{\widetilde{e_j}}X, \widetilde{e_j})g(Y, Z)\big]{\rm tr}[\texttt{id}],\nonumber\\
&\sum_{j=1}^{n}{\rm tr}[c(\widetilde{e_j})\nabla^{\bigwedge^*T^*M}_{\widetilde{e_j}}(A)]=\sum_{j=1}^{n}{\rm tr}[c(\widetilde{e_j})\nabla^{\bigwedge^*T^*M}_{\widetilde{e_j}}\big(c(X)c(Y)c(Z)\big)]\\
&=\big[\sum_{j=1}^{n}g(\widetilde{e_j}, Z)g(\nabla^{TM}_{\widetilde{e_j}}X,Y)-\sum_{j=1}^{n}g(\nabla^{TM}_{\widetilde{e_j}}X, Z)g(\widetilde{e_j}, Y)+\sum_{j=1}^{n}g(Y, Z)g(\widetilde{e_j}, \nabla^{TM}_{\widetilde{e_j}}X)\nonumber\\
&+\sum_{j=1}^{n}g(\widetilde{e_j}, Z)g(X,\nabla^{TM}_{\widetilde{e_j}}Y)-\sum_{j=1}^{n}g(X, Z)g(\widetilde{e_j}, \nabla^{TM}_{\widetilde{e_j}}Y)+\sum_{j=1}^{n}g(\nabla^{TM}_{\widetilde{e_j}}Y, Z)g(\widetilde{e_j}, X)\nonumber\\
&+\sum_{j=1}^{n}g(\widetilde{e_j}, \nabla^{TM}_{\widetilde{e_j}}Z)g(X, Y)-\sum_{j=1}^{n}g(X, \nabla^{TM}_{\widetilde{e_j}}Z)g(\widetilde{e_j}, Y)+\sum_{j=1}^{n}g(Y, \nabla^{TM}_{\widetilde{e_j}}Z)g(\widetilde{e_j}, X)\big]{\rm tr}[\texttt{id}],\nonumber
\end{align}
\begin{align}
&{\rm tr}[c(dx_n)A^*]=-g(\partial{x_n}, X)g(Y, Z){\rm tr}[\texttt{id}]+g(\partial{x_n}, Y)g(X, Z){\rm tr}[\texttt{id}]-g(\partial{x_n}, Z)g(X, Y){\rm tr}[\texttt{id}],\\
&{\rm tr}[c(dx_n)A]=g(\partial{x_n}, X)g(Y, Z){\rm tr}[\texttt{id}]-g(\partial{x_n}, Y)g(X, Z){\rm tr}[\texttt{id}]+g(\partial{x_n}, Z)g(X, Y){\rm tr}[\texttt{id}],
\end{align}
then we have
\begin{align}
&\int_M\int_{|\xi|=1}{\rm
trace}_{\wedge^*T^*M}[\sigma_{-4}(({{D}^*_{A}}{{D}_{A}})^{-1})]\sigma(\xi)dx=512\pi^2\int_{M}\bigg(-\frac{1}{12}s+2|X|^2|Y|^2|Z|^2\nonumber\\
&-\sum_{j=1}^{4}g(\nabla^{TM}_{\widetilde{e_j}}Z, \widetilde{e_j})g(X,Y)+\sum_{j=1}^{4}g(Y, \widetilde{e_j})g(X,\nabla^{TM}_{\widetilde{e_j}}Z)-\sum_{j=1}^{4}g(X, \widetilde{e_j})g(Y,\nabla^{TM}_{\widetilde{e_j}}Z)\nonumber\\
&-\sum_{j=1}^{4}g(Z, \widetilde{e_j})g(X,\nabla^{TM}_{\widetilde{e_j}}Y)+\sum_{j=1}^{4}g(\nabla^{TM}_{\widetilde{e_j}}Y, \widetilde{e_j})g(X,Z)-\sum_{j=1}^{4}g(X, \widetilde{e_j})g(\nabla^{TM}_{\widetilde{e_j}}Y,Z)\nonumber\\
&-\sum_{j=1}^{4}g(Z, \widetilde{e_j})g(\nabla^{TM}_{\widetilde{e_j}}X,Y)+\sum_{j=1}^{4}g(Y, \widetilde{e_j})g(\nabla^{TM}_{\widetilde{e_j}}X,Z)-\sum_{j=1}^{4}g(\nabla^{TM}_{\widetilde{e_j}}X, \widetilde{e_j})g(Y,Z)
\bigg)d{\rm Vol_{M}}\nonumber
\end{align}
and
\begin{align}
&\int_{\partial M}\Psi=\int_{\partial M} 8\pi\Omega_3 g(\partial{x_n}, X)g(Y,Z)d{\rm Vol_{\partial M}}-\int_{\partial M} 8\pi\Omega_3 g(\partial{x_n}, Y)g(X,Z)d{\rm Vol_{\partial M}}\\
&+\int_{\partial M} 8\pi\Omega_3 g(\partial{x_n}, Z)g(X,Y)d{\rm Vol_{\partial M}}.\nonumber
\end{align}
\begin{cor}
Let $M$ be a $4$-dimensional oriented
compact manifolds with the boundary $\partial M$ and the metric
$g^M$ as above, and let $A=c(X)c(Y)c(Z),$ then
\begin{align}
&\widetilde{{\rm Wres}}[\pi^+{{D}_{A}}^{-1}\circ\pi^+({D}^*_{A})^{-1}]=512\pi^2\int_{M}\bigg(-\frac{1}{12}s+2|X|^2|Y|^2|Z|^2\\
&-\sum_{j=1}^{4}g(\nabla^{TM}_{\widetilde{e_j}}Z, \widetilde{e_j})g(X,Y)+\sum_{j=1}^{4}g(Y, \widetilde{e_j})g(X,\nabla^{TM}_{\widetilde{e_j}}Z)-\sum_{j=1}^{4}g(X, \widetilde{e_j})g(Y,\nabla^{TM}_{\widetilde{e_j}}Z)\nonumber\\
&-\sum_{j=1}^{4}g(Z, \widetilde{e_j})g(X,\nabla^{TM}_{\widetilde{e_j}}Y)+\sum_{j=1}^{4}g(\nabla^{TM}_{\widetilde{e_j}}Y, \widetilde{e_j})g(X,Z)-\sum_{j=1}^{4}g(X, \widetilde{e_j})g(\nabla^{TM}_{\widetilde{e_j}}Y,Z)\nonumber\\
&-\sum_{j=1}^{4}g(Z, \widetilde{e_j})g(\nabla^{TM}_{\widetilde{e_j}}X,Y)+\sum_{j=1}^{4}g(Y, \widetilde{e_j})g(\nabla^{TM}_{\widetilde{e_j}}X,Z)-\sum_{j=1}^{4}g(\nabla^{TM}_{\widetilde{e_j}}X, \widetilde{e_j})g(Y,Z)
\bigg)d{\rm Vol_{M}}\nonumber\\
&+\int_{\partial M} 8\pi\Omega_3 g(\partial{x_n}, X)g(Y,Z)d{\rm Vol_{\partial M}}-\int_{\partial M} 8\pi\Omega_3 g(\partial{x_n}, Y)g(X,Z)d{\rm Vol_{\partial M}}\nonumber\\
&+\int_{\partial M} 8\pi\Omega_3 g(\partial{x_n}, Z)g(X,Y)d{\rm Vol_{\partial M}},\nonumber
\end{align}
where $s$ is the scalar curvature.
\end{cor}

When $A=\overline{c}(X)c(Y)c(Z),$ similar to Corollary 3.14, we have:
\begin{cor}
Let $M$ be a $4$-dimensional oriented
compact manifolds with the boundary $\partial M$ and the metric
$g^M$ as above, and let $A=\overline{c}(X)c(Y)c(Z),$ then
\begin{align}
&\widetilde{{\rm Wres}}[\pi^+{{D}_{A}}^{-1}\circ\pi^+({D}^*_{A})^{-1}]=512\pi^2\int_{M}\bigg(-\frac{1}{12}s+|X|^2|Y|^2|Z|^2-2|X|^2g(Y,Z)^2\bigg)d{\rm Vol_{M}},
\end{align}
where $s$ is the scalar curvature.
\end{cor}

When $A=\overline{c}(X)\overline{c}(Y)c(Z),$ we can get the following corollary:
\begin{cor}
Let $M$ be a $4$-dimensional oriented
compact manifolds with the boundary $\partial M$ and the metric
$g^M$ as above, and let $A=\overline{c}(X)\overline{c}(Y)c(Z),$ then
\begin{align}
&\widetilde{{\rm Wres}}[\pi^+{{D}_{A}}^{-1}\circ\pi^+({D}^*_{A})^{-1}]=512\pi^2\int_{M}\bigg(-\frac{1}{12}s+2|Z|^2g(X,Y)^2+\sum_{j=1}^{4}g(\nabla^{TM}_{\widetilde{e_j}}X,Y)g(Z,\widetilde{e_j})\\
&+\sum_{j=1}^{4}g(X,\nabla^{TM}_{\widetilde{e_j}}Y)g(Z,\widetilde{e_j})+\sum_{j=1}^{4}g(X,Y)g(\nabla^{TM}_{\widetilde{e_j}}Z,\widetilde{e_j})\bigg)d{\rm Vol_{M}}\nonumber\\
&-\int_{\partial M}8\pi\Omega_3 g(\partial{x_n},Z)g(X,Y)d{\rm Vol_{\partial M}},\nonumber
\end{align}
where $s$ is the scalar curvature.
\end{cor}

When $A=\overline{c}(X)\overline{c}(Y)\overline{c}(Z),$ we get:
\begin{cor}
Let $M$ be a $4$-dimensional oriented
compact manifolds with the boundary $\partial M$ and the metric
$g^M$ as above, and let $A=\overline{c}(X)\overline{c}(Y)\overline{c}(Z),$ then
\begin{align}
&\widetilde{{\rm Wres}}[\pi^+{{D}_{A}}^{-1}\circ\pi^+({D}^*_{A})^{-1}]=512\pi^2\int_{M}\bigg(-\frac{1}{12}s+3|X|^2|Y|^2|Z|^2-4|X|^2g(Y,Z)^2\\
&-4|Y|^2g(X,Z)^2-4|Z|^2g(X,Y)^2+8g(X,Y)g(X,Z)g(Y,Z)\bigg)d{\rm Vol_{M}},\nonumber
\end{align}
where $s$ is the scalar curvature.
\end{cor}

Next, we also prove the Kastler-Kalau-Walze type theorem for $4$-dimensional manifolds with boundary associated to ${{D}_{A}}^2$.
By (3.9) and (3.10), we will compute
\begin{equation}
\widetilde{{\rm Wres}}[\pi^+{{D}_{A}}^{-1}\circ\pi^+{{D}_{A}}^{-1}]=\int_M\int_{|\xi|=1}{\rm
trace}_{\wedge^*T^*M}[\sigma_{-4}({{D}_{A}}^{-2})]\sigma(\xi)dx+\int_{\partial M}\Phi,
\end{equation}
where
\begin{align}
\Phi &=\int_{|\xi'|=1}\int^{+\infty}_{-\infty}\sum^{\infty}_{j, k=0}\sum\frac{(-i)^{|\alpha|+j+k+1}}{\alpha!(j+k+1)!}
\times {\rm trace}_{\wedge^*T^*M}[\partial^j_{x_n}\partial^\alpha_{\xi'}\partial^k_{\xi_n}\sigma^+_{r}({{D}_{A}}^{-1})(x',0,\xi',\xi_n)\\
&\times\partial^\alpha_{x'}\partial^{j+1}_{\xi_n}\partial^k_{x_n}\sigma_{l}({{D}_{A}}^{-1})(x',0,\xi',\xi_n)]d\xi_n\sigma(\xi')dx',\nonumber
\end{align}
and the sum is taken over $r+l-k-j-|\alpha|=-3,~~r\leq -1,l\leq-1$.\\

By Theorem 2.2, we compute the interior of $\widetilde{{\rm Wres}}[\pi^+{{D}_{A}}^{-1}\circ\pi^+{{D}_{A}}^{-1}]$, then
\begin{align}
&\int_M\int_{|\xi|=1}{\rm
trace}_{\wedge^*T^*M}[\sigma_{-4}({{D}_{A}}^{-2})]\sigma(\xi)dx\\
&=32\pi^2\int_{M}{\rm tr}\bigg(-\frac{1}{12}s+A^2-\frac{1}{2}\sum_{j=1}^{n}Ac(\widetilde{e_j})Ac(\widetilde{e_j})\bigg)d{\rm Vol_{M}}.\nonumber
\end{align}

When $n=4$, then ${\rm tr}_{\wedge^*T^*M}[{\rm \texttt{id}}]={\rm dim}(\wedge^*(4))=16$, where ${\rm tr}$ as shorthand of ${\rm trace}$, the sum is taken over $
r+l-k-j-|\alpha|=-3,~~r\leq -1,l\leq-1,$ then we have the following five cases:
~\\

\noindent  {\bf case a)~I)}~$r=-1,~l=-1,~k=j=0,~|\alpha|=1$\\

\noindent By (3.94), we get
\begin{align}
\Phi_1=-\int_{|\xi'|=1}\int^{+\infty}_{-\infty}\sum_{|\alpha|=1}
 {\rm tr}[\partial^\alpha_{\xi'}\pi^+_{\xi_n}\sigma_{-1}({{D}_{A}}^{-1})\times
 \partial^\alpha_{x'}\partial_{\xi_n}\sigma_{-1}({{D}_{A}}^{-1})](x_0)d\xi_n\sigma(\xi')dx'.
\end{align}
\noindent  {\bf case a)~II)}~$r=-1,~l=-1,~k=|\alpha|=0,~j=1$\\

\noindent Likewise, we get
\begin{align}
\Phi_2=-\frac{1}{2}\int_{|\xi'|=1}\int^{+\infty}_{-\infty} {\rm
trace} [\partial_{x_n}\pi^+_{\xi_n}\sigma_{-1}({{D}_{A}}^{-1})\times
\partial_{\xi_n}^2\sigma_{-1}({{D}_{A}}^{-1})](x_0)d\xi_n\sigma(\xi')dx'.
\end{align}

\noindent  {\bf case a)~III)}~$r=-1,~l=-1,~j=|\alpha|=0,~k=1$\\

\noindent Observing (3.94), we get
\begin{align}
\Phi_3=-\frac{1}{2}\int_{|\xi'|=1}\int^{+\infty}_{-\infty}
{\rm trace} [\partial_{\xi_n}\pi^+_{\xi_n}\sigma_{-1}({{D}_{A}}^{-1})\times
\partial_{\xi_n}\partial_{x_n}\sigma_{-1}({{D}_{A}}^{-1})](x_0)d\xi_n\sigma(\xi')dx'.
\end{align}
By Lemma 3.7, we have $\sigma_{-1}({{D}_{A}}^{-1})=\sigma_{-1}(({D}^*_{A})^{-1})$.\\
In combination with the calculation,
\begin{align}
\Phi_1+\Phi_2+\Phi_3=0.
\end{align}

\noindent  {\bf case b)}~$r=-2,~l=-1,~k=j=|\alpha|=0$\\

\noindent By applying the formulae shown in (3.94), we get
\begin{align}
\Phi_4&=-i\int_{|\xi'|=1}\int^{+\infty}_{-\infty}{\rm trace} [\pi^+_{\xi_n}\sigma_{-2}({{D}_{A}}^{-1})\times
\partial_{\xi_n}\sigma_{-1}({{D}_{A}}^{-1})](x_0)d\xi_n\sigma(\xi')dx'.
\end{align}
By Lemma 3.7, we have $\sigma_{-1}({{D}_{A}}^{-1})=\sigma_{-1}(({D}^*_{A})^{-1})$.
Then, we have\\
\begin{align}
\Phi_4=\frac{9}{2}\pi h'(0)\Omega_3dx'+\frac{\pi}{4}\Omega_3{\rm tr}[Ac(dx_n)]dx',
\end{align}
where ${\rm \Omega_{4}}$ is the canonical volume of $S^{4}.$\\

\noindent {\bf  case c)}~$r=-1,~l=-2,~k=j=|\alpha|=0$\\

\noindent By (3.94), we can calculate
\begin{equation}
\Phi_5=-i\int_{|\xi'|=1}\int^{+\infty}_{-\infty}{\rm trace} [\pi^+_{\xi_n}\sigma_{-1}({{D}_{A}}^{-1})\times
\partial_{\xi_n}\sigma_{-2}({{D}_{A}}^{-1})](x_0)d\xi_n\sigma(\xi')dx'.
\end{equation}
By (3.5) and (3.6), we have
\begin{align}
\pi^+_{\xi_n}\sigma_{-1}({{D}_{A}}^{-1})(x_0)|_{|\xi'|=1}=\frac{c(\xi')+ic(dx_n)}{2(\xi_n-i)}.
\end{align}
Since
\begin{equation}
\sigma_{-2}({{D}_{A}}^{-1})(x_0)=\frac{c(\xi)\sigma_{0}({{D}_{A}})(x_0)c(\xi)}{|\xi|^4}+\frac{c(\xi)}{|\xi|^6}c(dx_n)\left[\partial_{x_n}[c(\xi')](x_0)|\xi|^2-c(\xi)h'(0)|\xi|^2_{\partial_M}\right],
\end{equation}
where
\begin{align}
\sigma_{0}({D}_{A})(x_0)&=
\frac{1}{4}\sum_{i,s,t}\omega_{s,t}(\widetilde{e_i})(x_{0})c(\widetilde{e_i})\overline{c}(\widetilde{e_s})\overline{c}(\widetilde{e_t})-\frac{1}{4}\sum_{i,s,t}\omega_{s,t}(\widetilde{e_i})(x_{0})c(\widetilde{e_i})c(\widetilde{e_s})c(\widetilde{e_t})+A\\
&=b_0^1(x_0)+b_0^2(x_0)+A,\nonumber
\end{align}
hence
\begin{align}
\partial_{\xi_n}\sigma_{-2}({{D}_{A}}^{-1})(x_0)|_{|\xi'|=1}&=
\partial_{\xi_n}\bigg(\frac{c(\xi)(b_0^1(x_0)+b_0^2(x_0)+A)c(\xi)}{|\xi|^4}\\
&+\frac{c(\xi)}{|\xi|^6}c(dx_n)\left[\partial_{x_n}[c(\xi')](x_0)|\xi|^2-c(\xi)h'(0)|\xi|^2_{\partial_M}\right]\bigg)\nonumber\\
&=\partial_{\xi_n}\bigg(\frac{c(\xi)b_0^2(x_0)c(\xi)}{|\xi|^4}+\frac{c(\xi)}{|\xi|^6}c(dx_n)\left[\partial_{x_n}[c(\xi')](x_0)|\xi|^2-c(\xi)h'(0)\right]\bigg)\nonumber\\
&+\partial_{\xi_n}\bigg(\frac{c(\xi)b_0^1(x_0)c(\xi)}{|\xi|^4}\bigg)+\partial_{\xi_n}\bigg(\frac{c(\xi)Ac(\xi)}{|\xi|^4}\bigg).\nonumber
\end{align}
By calculation, we have
\begin{align}
\partial_{\xi_n}\bigg(\frac{c(\xi)b_0^1(x_0)c(\xi)}{|\xi|^4}\bigg)
&=\frac{c(dx_n)b_0^1(x_0)c(\xi)}{|\xi|^4}+\frac{c(\xi)b_0^1(x_0)c(dx_n)}{|\xi|^4}-\frac{4\xi_nc(\xi)b_0^1(x_0)c(\xi)}{|\xi|^6};\\
\partial_{\xi_n}\bigg(\frac{c(\xi)Ac(\xi)}{|\xi|^4}\bigg)
&=\frac{c(dx_n)Ac(\xi)}{|\xi|^4}+\frac{c(\xi)Ac(dx_n)}{|\xi|^4}-\frac{4\xi_nc(\xi)Ac(\xi)}{|\xi|^6}.
\end{align}
For brevity, we denote
\begin{align}
q_{-2}^{1}=\frac{c(\xi)b_0^2(x_0)c(\xi)}{|\xi|^4}+\frac{c(\xi)}{|\xi|^6}c(dx_n)\left[\partial_{x_n}[c(\xi')](x_0)|\xi|^2-c(\xi)h'(0)\right],
\end{align}
then
\begin{align}
\partial_{\xi_n}(q_{-2}^{1})&=\frac{1}{(1+\xi_n^2)^3}[(2\xi_n-2\xi_n^3)c(dx_n)b_0^2(x_0)c(dx_n)+(1-3\xi_n^2)c(dx_n)b_0^2(x_0)c(\xi')\\
&+(1-3\xi_n^2)c(\xi')b_0^2(x_0)c(dx_n)-4\xi_nc(\xi')b_0^2(x_0)c(\xi')+(3\xi_n^2-1)\partial_{x_n}c(\xi')\nonumber\\
&-4\xi_nc(\xi')c(dx_n)\partial_{x_n}c(\xi')+2h'(0)c(\xi')+2h'(0)\xi_nc(dx_n)]\nonumber\\
&+6\xi_nh'(0)\frac{c(\xi)c(dx_n)c(\xi)}{(1+\xi^2_n)^4}.\nonumber
\end{align}
We calculate
\begin{align}
&{\rm tr}\bigg[\pi^+_{\xi_n}\sigma_{-1}({{D}_{A}}^{-1})(x_0)\times
\partial_{\xi_n}\bigg(\frac{c(\xi)b_0^1(x_0)c(\xi)}{|\xi|^4}\bigg)\bigg]|_{|\xi'|=1}\\
&=\frac{-1}{2(\xi_n-i)(\xi_n+i)^3}{\rm tr}[b_0^1(x_0)c(\xi')]+\frac{i}{2(\xi_n-i)(\xi_n+i)^3}{\rm tr}[b_0^1(x_0)c(dx_n)]\nonumber\\
&=\frac{-1}{2(\xi_n-i)(\xi_n+i)^3}{\rm tr}[b_0^1(x_0)c(\xi')],\nonumber
\end{align}
then
\begin{align}
&-i\int_{|\xi'|=1}\int^{+\infty}_{-\infty}{\rm trace}\bigg[\pi^+_{\xi_n}\sigma_{-1}({{D}_{A}}^{-1})\times
\partial_{\xi_n}\bigg(\frac{c(\xi)b_0^1c(\xi)}{|\xi|^4}\bigg)\bigg](x_0)d\xi_n\sigma(\xi')dx'\\
&=-i\int_{|\xi'|=1}\int^{+\infty}_{-\infty}\frac{-1}{2(\xi_n-i)(\xi_n+i)^3}{\rm tr}[b_0^1(x_0)c(\xi')]d\xi_n\sigma(\xi')dx'\nonumber\\
&=0.\nonumber
\end{align}
Likewise, we have
\begin{align}
&{\rm tr}\bigg[\pi^+_{\xi_n}\sigma_{-1}({{D}_{A}}^{-1})(x_0)\times
\partial_{\xi_n}\bigg(\frac{c(\xi)Ac(\xi)}{|\xi|^4}\bigg)\bigg]|_{|\xi'|=1}\\
&=\frac{-1}{2(\xi_n-i)(\xi_n+i)^3}{\rm tr}[Ac(\xi')]+\frac{i}{2(\xi_n-i)(\xi_n+i)^3}{\rm tr}[Ac(dx_n)],\nonumber
\end{align}
then by $\int_{|\xi'|=1}{\{\xi_{i_1}\cdot\cdot\cdot\xi_{i_{2d+1}}}\}\sigma(\xi')=0,$ we have
\begin{align}
&-i\int_{|\xi'|=1}\int^{+\infty}_{-\infty}{\rm trace}\bigg[\pi^+_{\xi_n}\sigma_{-1}({{D}_{A}}^{-1})\times
\partial_{\xi_n}\bigg(\frac{c(\xi)Ac(\xi)}{|\xi|^4}\bigg)\bigg](x_0)d\xi_n\sigma(\xi')dx'\\
&=-i\int_{|\xi'|=1}\int^{+\infty}_{-\infty}\frac{-1}{2(\xi_n-i)(\xi_n+i)^3}{\rm tr}[Ac(\xi')]d\xi_n\sigma(\xi')dx'\nonumber\\
&-i\int_{|\xi'|=1}\int^{+\infty}_{-\infty}\frac{i}{2(\xi_n-i)(\xi_n+i)^3}{\rm tr}[Ac(dx_n)]d\xi_n\sigma(\xi')dx'\nonumber\\
&=\frac{\Omega_3}{2}{\rm tr}[Ac(dx_n)]\int_{\Gamma^+}\frac{1}{(\xi_n-i)(\xi_n+i)^3}|_{\xi_n=i}dx'\nonumber\\
&=\frac{\Omega_3}{2}{\rm tr}[Ac(dx_n)]2\pi i[\frac{1}{(\xi_n+i)^3}]^{(1)}|_{\xi_n=i}dx'\nonumber\\
&=-\frac{\pi}{4}\Omega_3{\rm tr}[Ac(dx_n)]dx'.\nonumber
\end{align}
By (3.103) and (3.109), we have
\begin{align}
{\rm tr}[\pi^+_{\xi_n}\sigma_{-1}({{D}_{A}}^{-1})\times
\partial_{\xi_n}(q^1_{-2})](x_0)|_{|\xi'|=1}
=\frac{12h'(0)(i\xi^2_n+\xi_n-2i)}{(\xi-i)^3(\xi+i)^3}
+\frac{48h'(0)i\xi_n}{(\xi-i)^3(\xi+i)^4}.
\end{align}
We compute that
\begin{align}
&-i\int_{|\xi'|=1}\int^{+\infty}_{-\infty}{\rm tr}[\pi^+_{\xi_n}\sigma_{-1}({{D}_{A}}^{-1})\times
\partial_{\xi_n}(q^1_{-2})](x_0)d\xi_n\sigma(\xi')dx'\\
&=-i\int_{|\xi'|=1}\int^{+\infty}_{-\infty}\frac{12h'(0)(i\xi^2_n+\xi_n-2i)}{(\xi-i)^3(\xi+i)^3}+\frac{48h'(0)i\xi_n}{(\xi-i)^3(\xi+i)^4}d\xi_n\sigma(\xi')dx'\nonumber\\
&=-i\Omega_3\int_{\Gamma^+}\frac{12h'(0)(i\xi^2_n+\xi_n-2i)}{(\xi-i)^3(\xi+i)^3}+\frac{48h'(0)i\xi_n}{(\xi-i)^3(\xi+i)^4}d\xi_ndx'\nonumber\\
&=-i\Omega_3\frac{2\pi i}{2!}[\frac{12h'(0)(i\xi^2_n+\xi_n-2i)}{(\xi+i)^3}]^{(2)}|_{\xi_n=i}dx'-i\Omega_3\frac{2\pi i}{2!}[\frac{48h'(0)i\xi_n}{(\xi+i)^4}]^{(2)}|_{\xi_n=i}dx'\nonumber\\
&=-\frac{9}{2}\pi h'(0)\Omega_3dx'.\nonumber
\end{align}
Then,
\begin{align}
\Phi_5=-\frac{9}{2}\pi h'(0)\Omega_3dx'-\frac{\pi}{4}\Omega_3{\rm tr}[Ac(dx_n)]dx'.
\end{align}
So
\begin{align}
\Phi=\Phi_1+\Phi_2+\Phi_3+\Phi_4+\Phi_5=0.
\end{align}

\begin{thm}
Let $M$ be a $4$-dimensional oriented
compact manifold with the boundary $\partial M$ and the metric
$g^M$ as above, ${D}_{A}$ be the perturbation of the de Rham Hodge operator on $\widetilde{M}$, then
\begin{align}
&\widetilde{{\rm Wres}}[\pi^+{{D}_{A}}^{-1}\circ\pi^+{{D}_{A}}^{-1}]\nonumber\\
&=32\pi^2\int_{M}{\rm tr}\bigg(-\frac{1}{12}s+A^2-\frac{1}{2}\sum_{j=1}^{n}Ac(\widetilde{e_j})Ac(\widetilde{e_j})\bigg)d{\rm Vol_{M}}.\nonumber
\end{align}
where $s$ is the scalar curvature.
\end{thm}

We can directly state the following facts as a corollary of Theorem 3.18.

\begin{cor}
Let $M$ be a $4$-dimensional oriented
compact manifolds with the boundary $\partial M$ and the metric
$g^M$ as above, and let $A=c(X),$ then
\begin{align}
&\widetilde{{\rm Wres}}[\pi^+{{D}_{A}}^{-1}\circ\pi^+{{D}_{A}}^{-1}]=512\pi^2\int_{M}\bigg(-\frac{1}{12}s\bigg)d{\rm Vol_{M}},
\end{align}
where $s$ is the scalar curvature.
\end{cor}

When $A=\overline{c}(X),$ we can get the following corollary:
\begin{cor}
Let $M$ be a $4$-dimensional oriented
compact manifolds with the boundary $\partial M$ and the metric
$g^M$ as above, and let $A=\overline{c}(X),$ then
\begin{align}
\widetilde{{\rm Wres}}[\pi^+{{D}_{A}}^{-1}\circ\pi^+{{D}_{A}}^{-1}]&=512\pi^2\int_{M}\bigg(-\frac{1}{12}s-|X|^2\bigg)d{\rm Vol_{M}},
\end{align}
where $s$ is the scalar curvature.
\end{cor}

When $A=c(X)c(Y),$ similar to Corollary 3.20, we have :
\begin{cor}
Let $M$ be a $4$-dimensional oriented
compact manifolds with the boundary $\partial M$ and the metric
$g^M$ as above, and let $A=c(X)c(Y),$ then
\begin{align}
&\widetilde{{\rm Wres}}[\pi^+{{D}_{A}}^{-1}\circ\pi^+{{D}_{A}}^{-1}]=512\pi^2\int_{M}\bigg(-\frac{1}{12}s+|X|^2|Y|^2+2g(X,Y)^2\bigg)d{\rm Vol_{M}},
\end{align}
where $s$ is the scalar curvature.
\end{cor}

When $A=c(X)\overline{c}(Y),$ we compute $\widetilde{{\rm Wres}}[\pi^+{{D}_{A}}^{-1}\circ\pi^+{{D}_{A}}^{-1}].$
\begin{cor}
Let $M$ be a $4$-dimensional oriented
compact manifolds with the boundary $\partial M$ and the metric
$g^M$ as above, and let $A=c(X)\overline{c}(Y),$ then
\begin{align}
&\widetilde{{\rm Wres}}[\pi^+{{D}_{A}}^{-1}\circ\pi^+{{D}_{A}}^{-1}]=512\pi^2\int_{M}\bigg(-\frac{1}{12}s\bigg)d{\rm Vol_{M}},
\end{align}
where $s$ is the scalar curvature.
\end{cor}

A simple calculation shows that:
\begin{cor}
Let $M$ be a $4$-dimensional oriented
compact manifolds with the boundary $\partial M$ and the metric
$g^M$ as above, and let $A=\overline{c}(X)\overline{c}(Y),$ then
\begin{align}
&\widetilde{{\rm Wres}}[\pi^+{{D}_{A}}^{-1}\circ\pi^+{{D}_{A}}^{-1}]=512\pi^2\int_{M}\bigg(-\frac{1}{12}s-|X|^2|Y|^2+4g(X,Y)^2\bigg)d{\rm Vol_{M}},
\end{align}
where $s$ is the scalar curvature.
\end{cor}

When $A=c(X)c(Y)c(Z),$ we have:
\begin{cor}
Let $M$ be a $4$-dimensional oriented
compact manifolds with the boundary $\partial M$ and the metric
$g^M$ as above, and let $A=c(X)c(Y)c(Z),$ then
\begin{align}
&\widetilde{{\rm Wres}}[\pi^+{{D}_{A}}^{-1}\circ\pi^+{{D}_{A}}^{-1}]=512\pi^2\int_{M}\bigg(-\frac{1}{12}s\bigg)d{\rm Vol_{M}},
\end{align}
where $s$ is the scalar curvature.
\end{cor}

When $A=\overline{c}(X)c(Y)c(Z),$ we notice that:
\begin{cor}
Let $M$ be a $4$-dimensional oriented
compact manifolds with the boundary $\partial M$ and the metric
$g^M$ as above, and let $A=\overline{c}(X)c(Y)c(Z),$ then
\begin{align}
&\widetilde{{\rm Wres}}[\pi^+{{D}_{A}}^{-1}\circ\pi^+{{D}_{A}}^{-1}]=512\pi^2\int_{M}\bigg(-\frac{1}{12}s+|X|^2|Y|^2|Z|^2-2|X|^2g(Y,Z)^2\bigg)d{\rm Vol_{M}},
\end{align}
where $s$ is the scalar curvature.
\end{cor}

When $A=\overline{c}(X)\overline{c}(Y)c(Z),$ we have:
\begin{cor}
Let $M$ be a $4$-dimensional oriented
compact manifolds with the boundary $\partial M$ and the metric
$g^M$ as above, and let $A=\overline{c}(X)\overline{c}(Y)c(Z),$ then
\begin{align}
&\widetilde{{\rm Wres}}[\pi^+{{D}_{A}}^{-1}\circ\pi^+{{D}_{A}}^{-1}]=512\pi^2\int_{M}\bigg(-\frac{1}{12}s-2|X|^2|Y|^2|Z|^2+2|Z|^2g(X,Y)^2\bigg)d{\rm Vol_{M}},
\end{align}
where $s$ is the scalar curvature.
\end{cor}

When $A=\overline{c}(X)\overline{c}(Y)\overline{c}(Z),$ then  $A^*=\overline{c}(Z)\overline{c}(Y)\overline{c}(X).$\\
By computation, we have
\begin{align}
&{\rm tr}[A^2]={\rm tr}[|X|^2|Y|^2|Z|^2]=|X|^2|Y|^2|Z|^2{\rm tr}[\texttt{id}],\\
&{\rm tr}[c(dx_n)A^*]=0, {\rm tr}[c(dx_n)A]=0,\nonumber\\
&\sum_{j=1}^{n}{\rm tr}[Ac(\widetilde{e_j})Ac(\widetilde{e_j})]=\sum_{j=1}^{n}{\rm tr}[\overline{c}(X)\overline{c}(Y)\overline{c}(Z)\overline{c}(X)\overline{c}(Y)\overline{c}(Z)]\\
&=-\sum_{j=1}^{n}{\rm tr}[\overline{c}(Y)\overline{c}(X)\overline{c}(Z)\overline{c}(X)\overline{c}(Y)\overline{c}(Z)]+2\sum_{j=1}^{n}g(X, Y){\rm tr}[\overline{c}(Z)\overline{c}(X)\overline{c}(Y)\overline{c}(Z)]\nonumber\\
&=-n|X|^2|Y|^2|Z|^2{\rm tr}[\texttt{id}]+2n|X|^2g(Y, Z)^2{\rm tr}[\texttt{id}]+2n|Y|^2g(X, Z)^2{\rm tr}[\texttt{id}]\nonumber\\
&+2n|Z|^2g(X, Y)^2{\rm tr}[\texttt{id}]-4ng(X,Y)g(X, Z)(Y, Z){\rm tr}[\texttt{id}]\nonumber
\end{align}
We can calculate
\begin{align}
&\int_M\int_{|\xi|=1}{\rm
trace}_{\wedge^*T^*M}[\sigma_{-4}({{D}_{A}}^{-2})]\sigma(\xi)dx\\
&=32\pi^2\int_{M}{\rm tr}\bigg(-\frac{1}{12}s+A^2-\frac{1}{2}\sum_{j=1}^{n}Ac(\widetilde{e_j})Ac(\widetilde{e_j})\bigg)d{\rm Vol_{M}}\nonumber\\
&=32\pi^2\int_{M}\bigg(-\frac{1}{12}s+|X|^2|Y|^2|Z|^2+2|X|^2|Y|^2|Z|^2-4|X|^2g(Y, Z)^2-4|Y|^2g(X, Z)^2\nonumber\\
&-4|Z|^2g(X, Y)^2+8g(X,Y)g(X, Z)(Y, Z)\bigg){\rm tr}[\texttt{id}]d{\rm Vol_{M}}\nonumber\\
&=512\pi^2\int_{M}\bigg(-\frac{1}{12}s+3|X|^2|Y|^2|Z|^2-4|X|^2g(Y, Z)^2-4|Y|^2g(X, Z)^2\nonumber\\
&-4|Z|^2g(X, Y)^2+8g(X,Y)g(X, Z)(Y, Z)\bigg)d{\rm Vol_{M}},\nonumber
\end{align}
and
\begin{align}
\int_{\partial M} \Psi=0.
\end{align}
\begin{cor}
Let $M$ be a $4$-dimensional oriented
compact manifolds with the boundary $\partial M$ and the metric
$g^M$ as above, and let $A=\overline{c}(X)\overline{c}(Y)\overline{c}(Z),$ then
\begin{align}
&\widetilde{{\rm Wres}}[\pi^+{{D}_{A}}^{-1}\circ\pi^+{{D}_{A}}^{-1}]=512\pi^2\int_{M}\bigg(-\frac{1}{12}s+3|X|^2|Y|^2|Z|^2-4|X|^2g(Y,Z)^2\\
&-4|Y|^2g(X,Z)^2-4|Z|^2g(X,Y)^2+8g(X,Y)g(X,Z)g(Y,Z)\bigg)d{\rm Vol_{M}},\nonumber
\end{align}
where $s$ is the scalar curvature.
\end{cor}

\vskip 1 true cm

\section{A Kastler-Kalau-Walze type theorem for $6$-dimensional manifolds with boundary}
Firstly, we prove the Kastler-Kalau-Walze type theorems for $6$-dimensional manifolds with boundary. From \cite{Wa5}, we know that
\begin{equation}
\widetilde{{\rm Wres}}[\pi^+{{D}_{A}}^{-1}\circ\pi^+({D}^*_{A}{D}_{A}{D}^*_{A})^{-1}]=\int_M\int_{|\xi|=1}{\rm trace}_{\wedge^*T^*M}[\sigma_{-4}(({D}^*_{A}{D}_{A})^{-2})]\sigma(\xi)dx+\int_{{\partial}M}\overline{\Psi},
\end{equation}
where
\begin{align}
\overline{\Psi} &=\int_{|\xi'|=1}\int^{+\infty}_{-\infty}\sum^{\infty}_{j, k=0}\sum\frac{(-i)^{|\alpha|+j+k+1}}{\alpha!(j+k+1)!}
\times {\rm trace}_{\wedge ^*T^*M}[\partial^j_{x_n}\partial^\alpha_{\xi'}{\partial_t}^k_{\xi_n}\sigma^+_{r}({{D}_{A}}^{-1})(x',0,\xi',\xi_n)\\
&\times\partial^\alpha_{x'}\partial^{j+1}_{\xi_n}\partial^k_{x_n}\sigma_{l}(({D}^*_{A}{D}_{A}{D}^*_{A})^{-1})(x',0,\xi',\xi_n)]d\xi_n\sigma(\xi')dx',\nonumber
\end{align}
and the sum is taken over $r+\ell-k-j-|\alpha|-1=-6, \ r\leq-1, \ell\leq -3$.
By Theorem 2.2, we compute the interior term of (4.1), then
\begin{align}
&\int_M\int_{|\xi|=1}{\rm
trace}_{\wedge^*T^*M}[\sigma_{-4}(({D}^*_{A}{D}_{A})^{-2})]\sigma(\xi)dx\\
&=128\pi^3\int_{M}{\rm tr}\bigg(-\frac{1}{12}s+2A^*A-\frac{1}{4}\sum_{j=1}^{n}Ac(\widetilde{e_j})Ac(\widetilde{e_j})-\frac{1}{4}\sum_{j=1}^{n}A^*c(\widetilde{e_j})A^*c(\widetilde{e_j})\nonumber\\
&+\frac{1}{2}\sum_{j=1}^{n}\nabla^{\bigwedge^*T^*M}_{\widetilde{e_j}}(A^*)c(\widetilde{e_j})-\frac{1}{2}\sum_{j=1}^{n}c(\widetilde{e_j})\nabla^{\bigwedge^*T^*M}_{\widetilde{e_j}}(A)\bigg)d{\rm Vol_{M}}.\nonumber
\end{align}
Next, we compute $\int_{\partial M} \overline{\Psi}$. By computation, we get
\begin{align}
{D}^*_{A}{D}_{A}{D}^*_{A}
&=\sum^{n}_{i=1}c(\widetilde{e_i})\langle \widetilde{e_i},dx_{l}\rangle(-g^{ij}\partial_{l}\partial_{i}\partial_{j})
+\sum^{n}_{i=1}c(\widetilde{e_i})\langle \widetilde{e_i},dx_{l}\rangle \bigg\{-(\partial_{l}g^{ij})\partial_{i}\partial_{j}-g^{ij}\bigg(4(\sigma_{i}\\
&+a_{i})\partial_{j}-2\Gamma^{k}_{ij}\partial_{k}\bigg)\partial_{l}\bigg\}+\sum^{n}_{i=1}c(\widetilde{e_i})\langle \widetilde{e_i},dx_{l}\rangle \bigg\{-2(\partial_{l}g^{ij})(\sigma_{i}+a_{i})\partial_{j}+g^{ij}(\partial_{l}\Gamma^{k}_{ij})\partial_{k}\nonumber\\
&-2g^{ij}[(\partial_{l}\sigma_{i})+(\partial_{l}a_{i})]\partial_{j}+(\partial_{l}g^{ij})\Gamma^{k}_{ij}\partial_{k}+\sum_{j,k}\Big[\partial_{l}\Big(c(\widetilde{e_j})A+A^*c(\widetilde{e_j})\Big)\Big]\langle \widetilde{e_j},dx_{k}\rangle\partial_{k}\nonumber\\
&+\sum_{j,k}\Big(c(\widetilde{e_j})A+A^*c(\widetilde{e_j})\Big)\Big[\partial_{l}\langle \widetilde{e_j},dx_{k}\rangle\Big]\partial_{k}\bigg\}+\sum^{n}_{i=1}c(\widetilde{e_i})\langle c(\widetilde{e_i}),dx_{l}\rangle \partial_{l}\bigg\{-g^{ij}\Big[(\partial_{i}\sigma_{j})\nonumber\\
&+(\partial_{i}a_{j})+\sigma_{i}\sigma_{j}+\sigma_{i}a_{j}+a_{i}\sigma_{j}+a_{i}a_{j}-\Gamma_{ij}^{k}\sigma_{k}-\Gamma_{ij}^{k}a_{k}\Big]+\sum_{i,j}g^{ij}\Big[c(\partial_{i})\partial_{j}(A)\nonumber\\
&+c(\partial_{i})\sigma_{j}A+c(\partial_{i})a_{j}A+A^*c(\partial_{i})\sigma_{j}+A^*c(\partial_{i})a_{j}\Big]-\frac{1}{8}\sum_{i,j,k,l}R_{ijkl}\overline{c}(\widetilde{e_i})\overline{c}(\widetilde{e_j})c(\widetilde{e_k})c(\widetilde{e_l})\nonumber\\
&+\frac{1}{4}s+A^*A\bigg\}+\Big[\sum^{n}_{i=1}c(\widetilde{e_i})(\sigma_{i}+a_{i})+A^*\Big](-g^{ij}\partial_{i}\partial_{j})+\sum^{n}_{i=1}c(\widetilde{e_i})\langle \widetilde{e_i},dx_{l}\rangle \bigg\{2\sum_{j,k}\nonumber\\
&\Big[c(\widetilde{e_j})A+A^*c(\widetilde{e_j})\Big]\times\langle\widetilde{e_j},dx_{k}\rangle\bigg\}\partial_{l}\partial_{k}+\Big[\sum^{n}_{i=1}c(\widetilde{e_i})(\sigma_{i}+a_{i})+A^*\Big]\bigg\{-\sum_{i,j}g^{ij}\Big[2(\sigma_{i}\nonumber\\
&+a_{i})\partial_{j}-\Gamma_{ij}^{k}\partial_{k}+(\partial_{i}\sigma_{j})+(\partial_{i}a_{j})+\sigma_{i}\sigma_{j}+\sigma_{i}a_{j}+a_{i}\sigma_{j}+a_{i}a_{j}-\Gamma_{ij}^{k}\sigma_{k}-\Gamma_{ij}^{k}a_{k}\Big]\nonumber\\
&+\sum_{i,j}g^{ij}\Big(c(\partial_{i})A+A^*c(\partial_{i})\Big)\partial_{j}
+\sum_{i,j}g^{ij}\Big[c(\partial_{i})\partial_{j}(A)+c(\partial_{i})\sigma_{j}A+c(\partial_{i})a_{j}A\nonumber\\
&+A^*c(\partial_{i})\sigma_{j}+A^*c(\partial_{i})a_{j}\Big]-\frac{1}{8}\sum_{i,j,k,l}R_{ijkl}\overline{c}(\widetilde{e_i})\overline{c}(\widetilde{e_j})c(\widetilde{e_k})c(\widetilde{e_l})+\frac{1}{4}s+A^*A\bigg\}.\nonumber
\end{align}
Then, we obtain
\begin{lem}
The following identities hold:
\begin{align}
\sigma_2({D}^*_{A}{D}_{A}{D}^*_{A})
&=\sum_{i,j,l}c(dx_{l})\partial_{l}(g^{ij})\xi_{i}\xi_{j} +c(\xi)(4\sigma^k+4a^k-2\Gamma^k)\xi_{k}+2[|\xi|^2A-c(\xi)A^*c(\xi)]\\
&+\frac{1}{4}|\xi|^2\sum_{s,t,l}\omega_{s,t}(\widetilde{e_l})[c(\widetilde{e_l})\overline{c}(\widetilde{e_s})\overline{c}(\widetilde{e_t})-c(\widetilde{e_l})c(\widetilde{e_s})c(\widetilde{e_t})]+|\xi|^2A^*;\nonumber\\
\sigma_{3}({D}^*_{A}{D}_{A}{D}^*_{A})&=ic(\xi)|\xi|^{2}.\nonumber
\end{align}
\end{lem}
Write
\begin{align}
\sigma({D}^*_{A}{D}_{A}{D}^*_{A})&=p_3+p_2+p_1+p_0;
~\sigma(({D}^*_{A}{D}_{A}{D}^*_{A})^{-1})=\sum^{\infty}_{j=3}q_{-j}.
\end{align}
By the composition formula of pseudodifferential operators, we have
\begin{align}
1=\sigma(({D}^*_{A}{D}_{A}{D}^*_{A})\circ ({D}^*_{A}{D}_{A}{D}^*_{A})^{-1})&=
\sum_{\alpha}\frac{1}{\alpha!}\partial^{\alpha}_{\xi}
[\sigma({D}^*_{A}{D}_{A}{D}^*_{A})]D^{\alpha}_{x}
[({D}^*_{A}{D}_{A}{D}^*_{A})^{-1}] \\
&=(p_3+p_2+p_1+p_0)(q_{-3}+q_{-4}+q_{-5}+\cdots)\nonumber\\
&+\sum_j(\partial_{\xi_j}p_3+\partial_{\xi_j}p_2+\partial_{\xi_j}p_1+\partial_{\xi_j}p_0)\nonumber\\
&(D_{x_j}q_{-3}+D_{x_j}q_{-4}+D_{x_j}q_{-5}+\cdots) \nonumber\\
&=p_3q_{-3}+(p_3q_{-4}+p_2q_{-3}+\sum_j\partial_{\xi_j}p_3D_{x_j}q_{-3})+\cdots,\nonumber
\end{align}
by (4.7), we have

\begin{equation}
q_{-3}=p_3^{-1};~q_{-4}=-p_3^{-1}[p_2p_3^{-1}+\sum_j\partial_{\xi_j}p_3D_{x_j}(p_3^{-1})].
\end{equation}
By Lemma 4.1, we have some symbols of operators.
\begin{lem} The following identities hold:
\begin{align}
\sigma_{-3}(({D}^*_{A}{D}_{A}{D}^*_{A})^{-1})&=\frac{ic(\xi)}{|\xi|^{4}};\\
\sigma_{-4}(({D}^*_{A}{D}_{A}{D}^*_{A})^{-1})&=
\frac{c(\xi)\sigma_{2}({D}^*_{A}{D}_{A}{D}^*_{A})c(\xi)}{|\xi|^8}
+\frac{c(\xi)}{|\xi|^8}\Big(|\xi|^2c(dx_n)\partial_{x_n}c(\xi')
-2h'(0)c(dx_n)c(\xi)\nonumber\\
&+2\xi_{n}c(\xi)\partial_{x_n}c(\xi')+4\xi_{n}h'(0)\Big).\nonumber
\end{align}
\end{lem}
When $n=6$, then ${\rm tr}_{\wedge ^*T^*M}[\texttt{id}]=64$, where ${\rm tr}$ as shorthand of ${\rm trace}$.
Since the sum is taken over $r+\ell-k-j-|\alpha|-1=-6, \ r\leq-1, \ell\leq -3$, then we have the $\int_{\partial_{M}}\overline{\Psi}$
is the sum of the following five cases:\\
~\\
\noindent  {\bf case (a)~(I)}~$r=-1, l=-3, j=k=0, |\alpha|=1$.\\

\noindent By (4.2), we obtain
 \begin{equation}
\overline{\Psi}_1=-\int_{|\xi'|=1}\int^{+\infty}_{-\infty}\sum_{|\alpha|=1}{\rm trace}
[\partial^{\alpha}_{\xi'}\pi^{+}_{\xi_{n}}\sigma_{-1}({{D}_{A}}^{-1})
      \times\partial^{\alpha}_{x'}\partial_{\xi_{n}}\sigma_{-3}(({D}^*_{A}{D}_{A}{D}^*_{A})^{-1})](x_0)d\xi_n\sigma(\xi')dx'.
\end{equation}

\noindent  {\bf case (a)~(II)}~$r=-1, l=-3, |\alpha|=k=0, j=1$.\\

\noindent By (4.2), we have
  \begin{equation}
\overline{\Psi}_2=-\frac{1}{2}\int_{|\xi'|=1}\int^{+\infty}_{-\infty} {\rm
trace}[\partial_{x_{n}}\pi^{+}_{\xi_{n}}\sigma_{-1}({{D}_{A}}^{-1})\times\partial^{2}_{\xi_{n}}\sigma_{-3}(({D}^*_{A}{D}_{A}{D}^*_{A})^{-1})](x_0)d\xi_n\sigma(\xi')dx'.
\end{equation}

\noindent  {\bf case (a)~(III)}~$r=-1,l=-3,|\alpha|=j=0,k=1$.\\

\noindent It is easy to check that
 \begin{equation}
\overline{\Psi}_3=-\frac{1}{2}\int_{|\xi'|=1}\int^{+\infty}_{-\infty}{\rm trace}[\partial_{\xi_{n}}\pi^{+}_{\xi_{n}}\sigma_{-1}({{D}_{A}}^{-1})
      \times\partial_{\xi_{n}}\partial_{x_{n}}\sigma_{-3}(({D}^*_{A}{D}_{A}{D}^*_{A})^{-1})](x_0)d\xi_n\sigma(\xi')dx'.
\end{equation}

For \cite{Wa5}, we have
 \begin{equation}
 \overline{\Psi}_1+\overline{\Psi}_2+\overline{\Psi}_3=\frac{10}{2}\pi h'(0)\Omega_4dx'.
\end{equation}

\noindent  {\bf case (b)}~$r=-1,l=-4,|\alpha|=j=k=0$.\\

\noindent By observing (4.2), we have
 \begin{align}
\overline{\Psi}_4&=-i\int_{|\xi'|=1}\int^{+\infty}_{-\infty}{\rm trace}[\pi^{+}_{\xi_{n}}\sigma_{-1}({{D}_{A}}^{-1})
      \times\partial_{\xi_{n}}\sigma_{-4}(({D}^*_{A}{D}_{A}{D}^*_{A})^{-1})](x_0)d\xi_n\sigma(\xi')dx'\\
&=i\int_{|\xi'|=1}\int^{+\infty}_{-\infty}{\rm trace} [\partial_{\xi_n}\pi^+_{\xi_n}\sigma_{-1}({{D}_{A}}^{-1})\times
\sigma_{-4}(({D}^*_{A}{D}_{A}{D}^*_{A})^{-1})](x_0)d\xi_n\sigma(\xi')dx'.\nonumber
\end{align}

In the normal coordinate, $g^{ij}(x_{0})=\delta^{j}_{i}$ and $\partial_{x_{j}}(g^{\alpha\beta})(x_{0})=0$, if $j<n$; $\partial_{x_{j}}(g^{\alpha\beta})(x_{0})=h'(0)\delta^{\alpha}_{\beta}$, if $j=n$.
So by  \cite{Wa3}, when $k<n$, we have $\Gamma^{n}(x_{0})=\frac{5}{2}h'(0)$, $\Gamma^{k}(x_{0})=0$, $\delta^{n}(x_{0})=0$ and $\delta^{k}=\frac{1}{4}h'(0)c(\widetilde{e_k})c(\widetilde{e_n})$. Then, we obtain
\begin{align}
\sigma_{-4}(({D}^*_{A}{D}_{A}{D}^*_{A})^{-1})(x_{0})
&=\frac{1}{|\xi|^8}c(\xi)\Big(h'(0)c(\xi)\sum_{k<n}\xi_kc(\widetilde{e_k})c(\widetilde{e_n})-h'(0)c(\xi)\sum_{k<n}\xi_k\overline{c}(\widetilde{e_k})\overline{c}(\widetilde{e_n})\\
&-5h'(0)\xi_nc(\xi)+2[|\xi|^2A-c(\xi)A^*c(\xi)]+\frac{5}{4}|\xi|^2h'(0)c(\widetilde{e_i})\overline{c}(dx_n)\overline{c}(\widetilde{e_i})\nonumber\\
&-\frac{1}{4}|\xi|^2h'(0)c(dx_n)+|\xi|^2A^*\Big)c(\xi)+\frac{c(\xi)}{|\xi|^8}\Big(|\xi|^2c(dx_n)\partial_{x_n}c(\xi')\nonumber\\
&-2h'(0)c(dx_n)c(\xi)+2\xi_{n}c(\xi)\partial_{x_n}c(\xi')+4\xi_{n}h'(0)\Big),\nonumber
\end{align}
\begin{align}
\partial_{\xi_n}\pi^+_{\xi_n}\sigma_{-1}({{D}_{A}}^{-1})(x_0)|_{|\xi'|=1}=-\frac{c(\xi')+ic(dx_n)}{2(\xi_n-i)^2}.
\end{align}
By (4.15) and (4.16), we have
\begin{align}
&{\rm tr} [\partial_{\xi_n}\pi^+_{\xi_n}\sigma_{-1}({{D}_{A}}^{-1})\times
\sigma_{-4}({D}^*_{A}{D}_{A}{D}^*_{A})^{-1}](x_0)|_{|\xi'|=1}\\
&=-\frac{64h'(0)}{2(\xi_{n}-i)^{2}(1+\xi_{n}^{2})^{4}}\big(-\frac{15}{4}i+16\xi_{n}+\frac{19}{2}i\xi_{n}^{2}+8\xi_{n}^{3}+\frac{21}{4}i\xi_{n}^{4}\big)\nonumber\\
&+\frac{2+4i\xi_{n}-2\xi_{n}^{2}}{2(\xi_{n}-i)^{2}(1+\xi_{n}^{2})^{3}}{\rm tr}\big[Ac(\xi')\big]+\frac{-2i+4\xi_{n}+2i\xi_{n}^{2}}{2(\xi_{n}-i)^{2}(1+\xi_{n}^{2})^{3}}{\rm tr}\big[Ac(dx_n)\big]\nonumber\\
&+\frac{3+2i\xi_{n}+\xi_{n}^{2}}{2(\xi_{n}-i)^{2}(1+\xi_{n}^{2})^{3}}{\rm tr}\big[A^*c(\xi')\big]+\frac{i+2\xi_{n}+3i\xi_{n}^{2}}{2(\xi_{n}-i)^{2}(1+\xi_{n}^{2})^{3}}{\rm tr}\big[A^*c(dx_n)\big].\nonumber
\end{align}
Consequently,
\begin{align}
\overline{\Psi}_4&=
i\int_{|\xi'|=1}\int^{+\infty}_{-\infty}\frac{-64h'(0)(-\frac{15}{4}i+16\xi_{n}+\frac{19}{2}i\xi_{n}^{2}+8\xi_{n}^{3}+\frac{21}{4}i\xi_{n}^{4})}{2(\xi_{n}-i)^{2}(1+\xi_{n}^{2})^{4}}d\xi_n\sigma(\xi')dx'\\ &+i\int_{|\xi'|=1}\int^{+\infty}_{-\infty}\frac{-2i+4\xi_{n}+2i\xi_{n}^{2}}{2(\xi_{n}-i)^{2}(1+\xi_{n}^{2})^{3}}{\rm tr}\big[Ac(dx_n)\big]d\xi_n\sigma(\xi')dx'\nonumber\\
&+i\int_{|\xi'|=1}\int^{+\infty}_{-\infty}\frac{i+2\xi_{n}+3i\xi_{n}^{2}}{2(\xi_{n}-i)^{2}(1+\xi_{n}^{2})^{3}}{\rm tr}\big[A^*c(dx_n)\big]d\xi_n\sigma(\xi')dx'\nonumber\\
&=\frac{111}{2}\pi h'(0)\Omega_4dx'-\frac{3}{8}\pi\Omega_4{\rm tr}\big[Ac(dx_n)\big]dx'+\frac{1}{16}\pi\Omega_4{\rm tr}\big[A^*c(dx_n)\big]dx'.\nonumber
\end{align}

\noindent {\bf  case (c)}~$r=-2,l=-3,|\alpha|=j=k=0$.\\

\noindent By (4.2), we obtain
\begin{equation}
\overline{\Psi}_5=-i\int_{|\xi'|=1}\int^{+\infty}_{-\infty}{\rm trace}[\pi^{+}_{\xi_{n}}\sigma_{-2}({{D}_{A}}^{-1})
      \times\partial_{\xi_{n}}\sigma_{-3}(({D}^*_{A}{D}_{A}{D}^*_{A})^{-1})](x_0)d\xi_n\sigma(\xi')dx'.
\end{equation}
By Lemma 4.1 and Lemma 4.2, we have
\begin{align}
\sigma_{-2}({{D}_{A}}^{-1})(x_0)&=\frac{c(\xi)\sigma_{0}({{D}_{A}})(x_0)c(\xi)}{|\xi|^4}+\frac{c(\xi)}{|\xi|^6}\sum_jc(dx_j)\left[\partial_{x_j}(c(\xi))|\xi|^2-c(\xi)\partial_{x_j}(|\xi|^2)\right](x_0),
\end{align}
where
\begin{align}
\sigma_{0}({D}_{A})(x_0)&=
\frac{1}{4}\sum_{i,s,t}\omega_{s,t}(\widetilde{e_i})(x_{0})c(\widetilde{e_i})\overline{c}(\widetilde{e_s})\overline{c}(\widetilde{e_t})-\frac{1}{4}\sum_{i,s,t}\omega_{s,t}(\widetilde{e_i})(x_{0})c(\widetilde{e_i})c(\widetilde{e_s})c(\widetilde{e_t})+A.
\end{align}
On the other hand,
\begin{align}
\partial_{\xi_{n}}\sigma_{-3}(({D}^*_{A}{D}_{A}{D}^*_{A})^{-1})(x_{_{0}})|_{|\xi'|=1}=\frac{-8i\xi_{n}c(\xi')}{(1+\xi_{n}^{2})^{5}}+\frac{i(1-7\xi_{n}^{2})c(dx_{n})}{(1+\xi_{n}^{2})^{5}}.
\end{align}
It is easy to obtain that
\begin{align}
\pi^{+}_{\xi_{n}}\Big(\sigma_{-2}({{D}_{A}}^{-1})\Big)(x_{_{0}})|_{|\xi'|=1}
&=\pi^{+}_{\xi_{n}}\Big(\frac{c(\xi)\sigma_{0}({D}_{A})(x_{0})c(\xi)
+c(\xi)c(dx_{n})\partial_{x_{n}}[c(\xi')](x_{0})}{(1+\xi^{2}_{n})^{2}}\Big)\\
&-h'(0)\pi^{+}_{\xi_{n}}\Big(\frac{c(\xi)c(dx_{n})c(\xi)}{(1+\xi^{2}_{n})^{3}}\Big).\nonumber
\end{align}
We denote
 \begin{align}
\sigma_{0}({D}_{A})(x_0)=b_0^1(x_0)+b_0^2(x_0)+A.
\end{align}
Then, we obtain
\begin{align}
\pi^{+}_{\xi_{n}}\Big(\sigma_{-2}({{D}_{A}}^{-1})\Big)(x_{_{0}})|_{|\xi'|=1}
&=\pi^{+}_{\xi_{n}}\Big(\frac{c(\xi)b_0^2(x_0)c(\xi)
+c(\xi)c(dx_{n})\partial_{x_{n}}[c(\xi')](x_{0})}{(1+\xi^{2}_{n})^{2}}\Big)\\
&-h'(0)\pi^{+}_{\xi_{n}}\Big(\frac{c(\xi)c(dx_{n})c(\xi)}{(1+\xi^{2}_{n})^{3}}\Big)+\pi^+_{\xi_n}\Big(\frac{c(\xi)b_0^1(x_0)c(\xi)}{(1+\xi_n^2)^2}\Big)\nonumber\\
&+\pi^+_{\xi_n}\Big(\frac{c(\xi)Ac(\xi)}{(1+\xi_n^2)^2}\Big).\nonumber
\end{align}
Furthermore,
\begin{align}
\pi^+_{\xi_n}\Big(\frac{c(\xi)b_0^1(x_0)c(\xi)}{(1+\xi_n^2)^2}\Big)
&=\pi^+_{\xi_n}\Big(\frac{c(\xi')b_0^1(x_0)c(\xi')}{(1+\xi_n^2)^2}\Big)+\pi^+_{\xi_n}\Big( \frac{\xi_nc(\xi')b_0^1(x_0)c(dx_{n})}{(1+\xi_n^2)^2}\Big)\\
&+\pi^+_{\xi_n}\Big(\frac{\xi_nc(dx_{n})b_0^1(x_0)c(\xi')}{(1+\xi_n^2)^2}\Big)+\pi^+_{\xi_n}\Big(\frac{\xi_n^{2}c(dx_{n})b_0^1(x_0)c(dx_{n})}{(1+\xi_n^2)^2}\Big)\nonumber\\
&=-\frac{c(\xi')b_0^1(x_0)c(\xi')(2+i\xi_{n})}{4(\xi_{n}-i)^{2}}-\frac{ic(\xi')b_0^1(x_0)c(dx_{n})}{4(\xi_{n}-i)^{2}}\nonumber\\
&-\frac{ic(dx_{n})b_0^1(x_0)c(\xi')}{4(\xi_{n}-i)^{2}}-\frac{i\xi_{n}c(dx_{n})b_0^1(x_0)c(dx_{n})}{4(\xi_{n}-i)^{2}}.\nonumber
\end{align}
For the sake of convenience in writing,
\begin{align}
\pi^+_{\xi_n}\Big(\frac{c(\xi)b_0^2(x_0)c(\xi)+c(\xi)c(dx_n)\partial_{x_n}(c(\xi'))(x_0)}{(1+\xi_n^2)^2}\Big)-h'(0)\pi^+_{\xi_n}\Big(\frac{c(\xi)c(dx_n)c(\xi)}{(1+\xi_n)^3}\Big):= B_1-B_2,
\end{align}
where
\begin{align}
B_1&=\frac{-1}{4(\xi_n-i)^2}\big[(2+i\xi_n)c(\xi')b_0^2(x_0)c(\xi')+i\xi_nc(dx_n)b_0^2(x_0)c(dx_n) \\
&+(2+i\xi_n)c(\xi')c(dx_n)\partial_{x_n}c(\xi')+ic(dx_n)b_0^2(x_0)c(\xi')
+ic(\xi')b_0^2(x_0)c(dx_n)-i\partial_{x_n}c(\xi')\big]\nonumber
\end{align}
and
\begin{align}
B_2&=\frac{h'(0)}{2}\Big(\frac{c(dx_n)}{4i(\xi_n-i)}+\frac{c(dx_n)-ic(\xi')}{8(\xi_n-i)^2}
+\frac{3\xi_n-7i}{8(\xi_n-i)^3}\big(ic(\xi')-c(dx_n)\big)\Big).
\end{align}
By (4.22) and (4.28), we have
\begin{align}
&{\rm tr }[B_1\times\partial_{\xi_n}\sigma_{-3}(({D}^*_{A}{D}_{A}{D}^*_{A})^{-1})(x_0)]|_{|\xi'|=1}={\rm tr }\Big[ \frac{-1}{4(\xi_n-i)^2}\big[(2+i\xi_n)c(\xi')b_0^2(x_0)c(\xi')\\
&+i\xi_nc(dx_n)b_0^2(x_0)c(dx_n)+ic(dx_n)b_0^2(x_0)c(\xi')+(2+i\xi_n)c(\xi')c(dx_n)\partial_{x_n}c(\xi')\nonumber\\
&+ic(\xi')b_0^2(x_0)c(dx_n)-i\partial_{x_n}c(\xi')\big]\times \frac{-8i\xi_nc(\xi')+(i-7i\xi_n^{2})c(dx_n)}{(1+\xi_n^{2})^5}\Big] \nonumber\\
&=-32h'(0)\frac{-2i+9\xi_n+14i\xi_n^{2}-7\xi_n^{3}}{4(\xi_n-i)^2(1+\xi_n^2)^5}-64h'(0)\frac{-\frac{3}{2}i+12\xi_n+\frac{21}{2}i\xi_n^{2}}{4(\xi_n-i)^2(1+\xi_n^2)^5}\nonumber\\
&=8h'(0)\frac{-5-28i\xi_n+7\xi_n^{2}}{(\xi_n-i)^6(\xi_n+i)^5}.\nonumber
\end{align}
Similarly, we have
\begin{align}
&{\rm tr }[B_2\times\partial_{\xi_n}\sigma_{-3}(({D}^*_{A}{D}_{A}{D}^*_{A})^{-1})(x_0)]|_{|\xi'|=1}={\rm tr }\Big\{ \frac{h'(0)}{2}\Big[\frac{c(dx_n)}{4i(\xi_n-i)}+\frac{c(dx_n)-ic(\xi')}{8(\xi_n-i)^2}\\
&+\frac{3\xi_n-7i}{8(\xi_n-i)^3}[ic(\xi')-c(dx_n)]\Big]\times\frac{-8i\xi_nc(\xi')+(i-7i\xi_n^{2})c(dx_n)}{(1+\xi_n^{2})^5}\Big\} \nonumber\\
&=-8h'(0)\frac{4i-23\xi_n-14i\xi_n^{2}+7\xi_n^{3}}{(\xi_n-i)^7(\xi_n+i)^5},\nonumber
\end{align}
thus
\begin{align}
&{\rm tr }[B_1-B_2\times\partial_{\xi_n}\sigma_{-3}(({D}^*_{A}{D}_{A}{D}^*_{A})^{-1})(x_0)]|_{|\xi'|=1}
=8h'(0)\frac{9i-56\xi_n-49i\xi_n^{2}+14\xi_n^{3}}{(\xi_n-i)^7(\xi_n+i)^5}.
\end{align}
We calculate that
\begin{align}
&{\rm tr }\bigg[\pi^+_{\xi_n}\Big(\frac{c(\xi)b_0^1(x_0)c(\xi)}{(1+\xi_n^2)^2}\Big)\times
\partial_{\xi_n}\sigma_{-3}(({D}^*_{A}{D}_{A}{D}^*_{A})^{-1})(x_0)\bigg]\bigg|_{|\xi'|=1}\\
&=\frac{-2-16i\xi_n+14\xi_n^{2}}{4(\xi_n-i)^7(\xi_n+i)^5}{\rm tr}[b_0^1(x_0)c(\xi')].\nonumber
\end{align}
Similar calculations to (4.33), it is shown that
\begin{align}
&{\rm tr }\bigg[\pi^+_{\xi_n}\Big(\frac{c(\xi)Ac(\xi)}{(1+\xi_n^2)^2}\Big)\times
\partial_{\xi_n}\sigma_{-3}(({D}^*_{A}{D}_{A}{D}^*_{A})^{-1})(x_0)\bigg]\bigg|_{|\xi'|=1}\\
&=\frac{-2-16i\xi_n+14\xi_n^{2}}{4(\xi_n-i)^7(\xi_n+i)^5}{\rm tr}[Ac(\xi')]+\frac{-2i+16\xi_n+14i\xi_n^{2}}{4(\xi_n-i)^7(\xi_n+i)^5}{\rm tr}[Ac(dx_n)].\nonumber
\end{align}
By $\int_{|\xi'|=1}\xi_{1}\cdot\cdot\cdot\xi_{2q+1}\sigma(\xi')=0,$ we have\\
\begin{align}
\overline{\Psi}_5
&=-i\int_{|\xi'|=1}\int^{+\infty}_{-\infty}8h'(0)\frac{9i-56\xi_n-49i\xi_n^{2}+14\xi_n^{3}}{(\xi_n-i)^7(\xi_n+i)^5}d\xi_n\sigma(\xi')dx'\\
&-i\int_{|\xi'|=1}\int^{+\infty}_{-\infty}\frac{-2-16i\xi_n+14\xi_n^{2}}{4(\xi_n-i)^7(\xi_n+i)^5}{\rm tr}[b_0^1(x_0)c(\xi')]d\xi_n\sigma(\xi')dx'\nonumber\\
&-i\int_{|\xi'|=1}\int^{+\infty}_{-\infty}\frac{-2-16i\xi_n+14\xi_n^{2}}{4(\xi_n-i)^7(\xi_n+i)^5}{\rm tr}[Ac(\xi')]d\xi_n\sigma(\xi')dx'\nonumber\\
&-i\int_{|\xi'|=1}\int^{+\infty}_{-\infty}\frac{-2i+16\xi_n+14i\xi_n^{2}}{4(\xi_n-i)^7(\xi_n+i)^5}{\rm tr}[Ac(dx_n)]d\xi_n\sigma(\xi')dx'\nonumber\\
&=-i8h'(0)\Omega_4\frac{2\pi i}{6!}\Big[\frac{9i-56\xi_n-49i\xi_n^{2}+14\xi_n^{3}}{(\xi_n+i)^5}\Big]^{(5)}|_{\xi_n=i}dx'\nonumber\\
&-i\Omega_4{\rm tr}[Ac(dx_n)]\frac{2\pi i}{6!}\Big[\frac{-2i+16\xi_n+14i\xi_n^{2}}{4(\xi_n+i)^5}\Big]^{(5)}|_{\xi_n=i}dx'\nonumber\\
&=-\frac{105}{4}\pi h'(0)\Omega_4dx'+\frac{161}{512}\pi\Omega_4{\rm tr}[Ac(dx_n)]dx'.\nonumber
\end{align}

Now $\overline{\Psi}$ is the sum of the cases (a), (b) and (c), then
\begin{equation}
\overline{\Psi}=\frac{137}{4}\pi h'(0)\Omega_4dx'-\frac{31}{512}\pi\Omega_4{\rm tr}\big[Ac(dx_n)\big]dx'+\frac{1}{16}\pi\Omega_4{\rm tr}\big[A^*c(dx_n)\big]dx'.
\end{equation}

\begin{thm}
Let $M$ be a $6$-dimensional
compact oriented manifold with the boundary $\partial M$ and the metric
$g^M$ as above, ${D}_{A}$ be the perturbation of the de Rham Hodge operator on $\widetilde{M}$, then
\begin{align}
&\widetilde{{\rm Wres}}[\pi^+{{D}_{A}}^{-1}\circ\pi^+({D}^*_{A}{D}_{A}{D}^*_{A})^{-1}]=128\pi^3\int_{M}{\rm tr}\bigg(-\frac{1}{12}s+2A^*A-\frac{1}{4}\sum_{j=1}^{n}Ac(\widetilde{e_j})Ac(\widetilde{e_j})\\
&-\frac{1}{4}\sum_{j=1}^{n}A^*c(\widetilde{e_j})A^*c(\widetilde{e_j})+\frac{1}{2}\sum_{j=1}^{n}\nabla^{\bigwedge^*T^*M}_{\widetilde{e_j}}(A^*)c(\widetilde{e_j})-\frac{1}{2}\sum_{j=1}^{n}c(\widetilde{e_j})\nabla^{\bigwedge^*T^*M}_{\widetilde{e_j}}(A)\bigg)d{\rm Vol_{M}}\nonumber\\
&+\int_{\partial M}\frac{137}{4}\pi h'(0)\Omega_4d{\rm Vol_{\partial M}}-\int_{\partial M}\frac{31}{512}\pi\Omega_4{\rm tr}\big[Ac(dx_n)\big]d{\rm Vol_{\partial M}}\nonumber\\
&+\int_{\partial M}\frac{1}{16}\pi\Omega_4{\rm tr}\big[A^*c(dx_n)\big]d{\rm Vol_{\partial M}},\nonumber
\end{align}
where $s$ is the scalar curvature.
\end{thm}

We can state the following facts as a corollary of Theorem 4.3.
\begin{cor}
Let $M$ be a $6$-dimensional oriented
compact manifolds with the boundary $\partial M$ and the metric
$g^M$ as above, and let $A=c(X),$ then
\begin{align}
&\widetilde{{\rm Wres}}[\pi^+{{D}_{A}}^{-1}\circ\pi^+({D}^*_{A}{D}_{A}{D}^*_{A})^{-1}]=8192\pi^3\int_{M}\bigg(-\frac{1}{12}s+4|X|^2+\sum_{j=1}^{n}g(\nabla^{TM}_{\widetilde{e_j}}X,\widetilde{e_j})\bigg)d{\rm Vol_{M}}\\
&+\int_{\partial M}\frac{137}{4}\pi h'(0)\Omega_4d{\rm Vol_{\partial M}}+\int_{\partial M}\frac{63}{8}\pi\Omega_4g(\partial{x_n}, X)d{\rm Vol_{\partial M}},\nonumber
\end{align}
where $s$ is the scalar curvature.
\end{cor}

When $A=\overline{c}(X),$ we can get the following corollary:
\begin{cor}
Let $M$ be a $6$-dimensional oriented
compact manifolds with the boundary $\partial M$ and the metric
$g^M$ as above, and let $A=\overline{c}(X),$ then
\begin{align}
&\widetilde{{\rm Wres}}[\pi^+{{D}_{A}}^{-1}\circ\pi^+({D}^*_{A}{D}_{A}{D}^*_{A})^{-1}]=8192\pi^3\int_{M}\bigg(-\frac{1}{12}s-|X|^2\bigg)d{\rm Vol_{M}}\\
&+\int_{\partial M}\frac{137}{4}\pi h'(0)\Omega_4d{\rm Vol_{\partial M}},\nonumber
\end{align}
where $s$ is the scalar curvature.
\end{cor}

When $A=c(X)c(Y),$ we compute $\widetilde{{\rm Wres}}[\pi^+{{D}_{A}}^{-1}\circ\pi^+({D}^*_{A}{D}_{A}{D}^*_{A})^{-1}].$
\begin{cor}
Let $M$ be a $6$-dimensional oriented
compact manifolds with the boundary $\partial M$ and the metric
$g^M$ as above, and let $A=c(X)c(Y),$ then
\begin{align}
&\widetilde{{\rm Wres}}[\pi^+{{D}_{A}}^{-1}\circ\pi^+({D}^*_{A}{D}_{A}{D}^*_{A})^{-1}]=8192\pi^3\int_{M}\bigg(-\frac{1}{12}s+|X|^2|Y|^2+4g(X,Y)^2\bigg)d{\rm Vol_{M}}\nonumber\\
&+\int_{\partial M}\frac{137}{4}\pi h'(0)\Omega_4d{\rm Vol_{\partial M}},\nonumber
\end{align}
where $s$ is the scalar curvature.
\end{cor}

When $A=c(X)\overline{c}(Y),$ we obtain:
\begin{cor}
Let $M$ be a $6$-dimensional oriented
compact manifolds with the boundary $\partial M$ and the metric
$g^M$ as above, and let $A=c(X)\overline{c}(Y),$ then
\begin{align}
&\widetilde{{\rm Wres}}[\pi^+{{D}_{A}}^{-1}\circ\pi^+({D}^*_{A}{D}_{A}{D}^*_{A})^{-1}]=8192\pi^3\int_{M}\bigg(-\frac{1}{12}s+4|X|^2|Y|^2\bigg)d{\rm Vol_{M}}\\
&+\int_{\partial M}\frac{137}{4}\pi h'(0)\Omega_4d{\rm Vol_{\partial M}},\nonumber
\end{align}
where $s$ is the scalar curvature.
\end{cor}

When $A=\overline{c}(X)\overline{c}(Y),$ similar to Corollary 4.7, we can get:
\begin{cor}
Let $M$ be a $6$-dimensional oriented
compact manifolds with the boundary $\partial M$ and the metric
$g^M$ as above, and let $A=\overline{c}(X)\overline{c}(Y),$ then
\begin{align}
&\widetilde{{\rm Wres}}[\pi^+{{D}_{A}}^{-1}\circ\pi^+({D}^*_{A}{D}_{A}{D}^*_{A})^{-1}]=8192\pi^3\int_{M}\bigg(-\frac{1}{12}s-|X|^2|Y|^2+6g(X,Y)^2\bigg)d{\rm Vol_{M}}\\
&+\int_{\partial M}\frac{137}{4}\pi h'(0)\Omega_4d{\rm Vol_{\partial M}},\nonumber
\end{align}
where $s$ is the scalar curvature.
\end{cor}

When $A=c(X)c(Y)c(Z),$ we compute $\widetilde{{\rm Wres}}[\pi^+{{D}_{A}}^{-1}\circ\pi^+({D}^*_{A}{D}_{A}{D}^*_{A})^{-1}].$
\begin{cor}
Let $M$ be a $6$-dimensional oriented
compact manifolds with the boundary $\partial M$ and the metric
$g^M$ as above, and let $A=c(X)c(Y)c(Z),$ then
\begin{align}
&\widetilde{{\rm Wres}}[\pi^+{{D}_{A}}^{-1}\circ\pi^+({D}^*_{A}{D}_{A}{D}^*_{A})^{-1}]=8192\pi^3\int_{M}\bigg(-\frac{1}{12}s+2|X|^2|Y|^2|Z|^2\\
&+2|X|^2g(Y,Z)^2+2|Y|^2g(X,Z)^2+2|Z|^2g(X,Y)^2-4g(X,Y)g(X,Z)g(Y,Z)\nonumber\\
&-\sum_{j=1}^{6}g(\nabla^{TM}_{\widetilde{e_j}}Z, \widetilde{e_j})g(X,Y)+\sum_{j=1}^{6}g(Y, \widetilde{e_j})g(X,\nabla^{TM}_{\widetilde{e_j}}Z)-\sum_{j=1}^{6}g(X, \widetilde{e_j})g(Y,\nabla^{TM}_{\widetilde{e_j}}Z)\nonumber\\
&-\sum_{j=1}^{6}g(Z, \widetilde{e_j})g(X,\nabla^{TM}_{\widetilde{e_j}}Y)+\sum_{j=1}^{6}g(\nabla^{TM}_{\widetilde{e_j}}Y, \widetilde{e_j})g(X,Z)-\sum_{j=1}^{6}g(X, \widetilde{e_j})g(\nabla^{TM}_{\widetilde{e_j}}Y,Z)\nonumber\\
&-\sum_{j=1}^{6}g(Z, \widetilde{e_j})g(\nabla^{TM}_{\widetilde{e_j}}X,Y)+\sum_{j=1}^{6}g(Y, \widetilde{e_j})g(\nabla^{TM}_{\widetilde{e_j}}X,Z)-\sum_{j=1}^{6}g(\nabla^{TM}_{\widetilde{e_j}}X, \widetilde{e_j})g(Y,Z)
\bigg)d{\rm Vol_{M}}\nonumber\\
&+\int_{\partial M}\frac{137}{4}\pi h'(0)\Omega_4d{\rm Vol_{\partial M}}+\int_{\partial M}\frac{1}{8}\pi\Omega_4g(\partial{x_n}, X)g(Y,Z)d{\rm Vol_{\partial M}}\nonumber\\
&-\int_{\partial M} \frac{1}{8}\pi\Omega_4g(\partial{x_n}, Y)g(X,Z)d{\rm Vol_{\partial M}}+\int_{\partial M} \frac{1}{8}\pi\Omega_4 g(\partial{x_n}, Z)g(X,Y)d{\rm Vol_{\partial M}},\nonumber
\end{align}
where $s$ is the scalar curvature.
\end{cor}

When $A=\overline{c}(X)c(Y)c(Z),$ we get:
\begin{cor}
Let $M$ be a $6$-dimensional oriented
compact manifolds with the boundary $\partial M$ and the metric
$g^M$ as above, and let $A=\overline{c}(X)c(Y)c(Z),$ then
\begin{align}
&\widetilde{{\rm Wres}}[\pi^+{{D}_{A}}^{-1}\circ\pi^+({D}^*_{A}{D}_{A}{D}^*_{A})^{-1}]\\
&=8192\pi^3\int_{M}\bigg(-\frac{1}{12}s+3|X|^2|Y|^2|Z|^2-4|X|^2g(Y,Z)^2\bigg)d{\rm Vol_{M}}+\int_{\partial M}\frac{137}{4}\pi h'(0)\Omega_4d{\rm Vol_{\partial M}},\nonumber
\end{align}
where $s$ is the scalar curvature.
\end{cor}

When $A=\overline{c}(X)\overline{c}(Y)c(Z),$ we conclude that:
\begin{cor}
Let $M$ be a $6$-dimensional oriented
compact manifolds with the boundary $\partial M$ and the metric
$g^M$ as above, and let $A=\overline{c}(X)\overline{c}(Y)c(Z),$ then
\begin{align}
&\widetilde{{\rm Wres}}[\pi^+{{D}_{A}}^{-1}\circ\pi^+({D}^*_{A}{D}_{A}{D}^*_{A})^{-1}]=8192\pi^3\int_{M}\bigg(-\frac{1}{12}s+4|Z|^2g(X,Y)^2\\
&+\sum_{j=1}^{6}g(\nabla^{TM}_{\widetilde{e_j}}X,Y)g(Z,\widetilde{e_j})+\sum_{j=1}^{6}g(X,\nabla^{TM}_{\widetilde{e_j}}Y)g(Z,\widetilde{e_j})+\sum_{j=1}^{6}g(X,Y)g(\nabla^{TM}_{\widetilde{e_j}}Z,\widetilde{e_j})\bigg)d{\rm Vol_{M}}\nonumber\\
&+\int_{\partial M}\frac{137}{4}\pi h'(0)\Omega_4d{\rm Vol_{\partial M}}+\int_{\partial M}\frac{63}{8}\pi\Omega_4g(\partial{x_n}, Z)g(X, Y)d{\rm Vol_{\partial M}},\nonumber
\end{align}
where $s$ is the scalar curvature.
\end{cor}

When $A=\overline{c}(X)\overline{c}(Y)\overline{c}(Z),$ we find that:
\begin{cor}
Let $M$ be a $6$-dimensional oriented
compact manifolds with the boundary $\partial M$ and the metric
$g^M$ as above, and let $A=\overline{c}(X)\overline{c}(Y)\overline{c}(Z),$ then
\begin{align}
&\widetilde{{\rm Wres}}[\pi^+{{D}_{A}}^{-1}\circ\pi^+({D}^*_{A}{D}_{A}{D}^*_{A})^{-1}]=8192\pi^3\int_{M}\bigg(-\frac{1}{12}s+5|X|^2|Y|^2|Z|^2-6|X|^2g(Y,Z)^2\\
&-6|Y|^2g(X,Z)^2-6|Z|^2g(X,Y)^2+12g(X,Y)g(X,Z)g(Y,Z)\bigg)d{\rm Vol_{M}}\nonumber\\
&+\int_{\partial M}\frac{137}{4}\pi h'(0)\Omega_4d{\rm Vol_{\partial M}},\nonumber
\end{align}
where $s$ is the scalar curvature.
\end{cor}

Next, we prove the Kastler-Kalau-Walze type theorem for $6$-dimensional manifold with boundary associated to ${{D}_{A}}^{3}$. From \cite{Wa5}, we know that

\begin{equation}
\widetilde{{\rm Wres}}[\pi^+{{D}_{A}}^{-1}\circ\pi^+{{D}_{A}}^{-3}]=\int_M\int_{|\xi|=1}{\rm
trace}_{\wedge^*T^*M}[\sigma_{-4}({{D}_{A}}^{-4})]\sigma(\xi)dx+\int_{\partial M}\overline{\Phi},
\end{equation}
where $\widetilde{{\rm Wres}}$ denote noncommutative residue on minifolds with boundary,
\begin{align}
\overline{\Phi} &=\int_{|\xi'|=1}\int^{+\infty}_{-\infty}\sum^{\infty}_{j, k=0}\sum\frac{(-i)^{|\alpha|+j+k+1}}{\alpha!(j+k+1)!}
\times {\rm trace}_{{\wedge^*T^*M}}[\partial^j_{x_n}\partial^\alpha_{\xi'}\partial^k_{\xi_n}\sigma^+_{r}({{D}_{A}}^{-1})(x',0,\xi',\xi_n)\\
&\times\partial^\alpha_{x'}\partial^{j+1}_{\xi_n}\partial^k_{x_n}\sigma_{l}
({{D}_{A}}^{-3})(x',0,\xi',\xi_n)]d\xi_n\sigma(\xi')dx',\nonumber
\end{align}
and the sum is taken over $r+\ell-k-j-|\alpha|-1=-6, \ r\leq-1, \ell\leq -3$.

By Theorem 2.2, we compute the interior term of (4.46)
\begin{align}
&\int_M\int_{|\xi|=1}{\rm trace}_{\wedge^*T^*M}[\sigma_{-4}({{D}_{A}}^{-4})]\sigma(\xi)dx\\
&=128\pi^3\int_{M}{\rm tr}\bigg(-\frac{1}{12}s+2A^2-\frac{1}{2}\sum_{j=1}^{n}Ac(\widetilde{e_j})Ac(\widetilde{e_j})\bigg)d{\rm Vol_{M}}.\nonumber
\end{align}

So we only need to compute $\int_{\partial M} \overline{\Phi}$. Let us now turn to compute the specification of
${{D}_{A}}^3$.
\begin{align}
{{D}_{A}}^3
&=\sum^{n}_{i=1}c(\widetilde{e_i})\langle \widetilde{e_i},dx_{l}\rangle(-g^{ij}\partial_{l}\partial_{i}\partial_{j})
+\sum^{n}_{i=1}c(\widetilde{e_i})\langle \widetilde{e_i},dx_{l}\rangle \bigg\{-(\partial_{l}g^{ij})\partial_{i}\partial_{j}-g^{ij}\bigg(4(\sigma_{i}\\
&+a_{i})\partial_{j}-2\Gamma^{k}_{ij}\partial_{k}\bigg)\partial_{l}\bigg\}+\sum^{n}_{i=1}c(\widetilde{e_i})\langle \widetilde{e_i},dx_{l}\rangle \bigg\{-2(\partial_{l}g^{ij})(\sigma_{i}+a_{i})\partial_{j}+g^{ij}(\partial_{l}\Gamma^{k}_{ij})\partial_{k}\nonumber\\
&-2g^{ij}[(\partial_{l}\sigma_{i})+(\partial_{l}a_{i})]\partial_{j}+(\partial_{l}g^{ij})\Gamma^{k}_{ij}\partial_{k}+\sum_{j,k}\Big[\partial_{l}\Big(c(\widetilde{e_j})A+Ac(\widetilde{e_j})\Big)\Big]\langle \widetilde{e_j},dx_{k}\rangle\partial_{k}\nonumber\\
&+\sum_{j,k}\Big(c(\widetilde{e_j})A+Ac(\widetilde{e_j})\Big)\Big[\partial_{l}\langle \widetilde{e_j},dx_{k}\rangle\Big]\partial_{k}\bigg\}+\sum^{n}_{i=1}c(\widetilde{e_i})\langle c(\widetilde{e_i}),dx_{l}\rangle \partial_{l}\bigg\{-g^{ij}\Big[(\partial_{i}\sigma_{j})\nonumber\\
&+(\partial_{i}a_{j})+\sigma_{i}\sigma_{j}+\sigma_{i}a_{j}+a_{i}\sigma_{j}+a_{i}a_{j}-\Gamma_{ij}^{k}\sigma_{k}-\Gamma_{ij}^{k}a_{k}\Big]+\sum_{i,j}g^{ij}\Big[c(\partial_{i})\partial_{j}(A)\nonumber\\
&+c(\partial_{i})\sigma_{j}A+c(\partial_{i})a_{j}A+Ac(\partial_{i})\sigma_{j}+Ac(\partial_{i})a_{j}\Big]-\frac{1}{8}\sum_{i,j,k,l}R_{ijkl}\overline{c}(\widetilde{e_i})\overline{c}(\widetilde{e_j})c(\widetilde{e_k})c(\widetilde{e_l})\nonumber\\
&+\frac{1}{4}s+A^2\bigg\}+\Big[\sum^{n}_{i=1}c(\widetilde{e_i})(\sigma_{i}+a_{i})+A\Big](-g^{ij}\partial_{i}\partial_{j})+\sum^{n}_{i=1}c(\widetilde{e_i})\langle \widetilde{e_i},dx_{l}\rangle \bigg\{2\sum_{j,k}\nonumber\\
&\Big[c(\widetilde{e_j})A+Ac(\widetilde{e_j})\Big]\times\langle\widetilde{e_j},dx_{k}\rangle\bigg\}\partial_{l}\partial_{k}+\Big[\sum^{n}_{i=1}c(\widetilde{e_i})(\sigma_{i}+a_{i})+A\Big]\bigg\{-\sum_{i,j}g^{ij}\Big[2(\sigma_{i}\nonumber\\
&+a_{i})\partial_{j}-\Gamma_{ij}^{k}\partial_{k}+(\partial_{i}\sigma_{j})+(\partial_{i}a_{j})+\sigma_{i}\sigma_{j}+\sigma_{i}a_{j}+a_{i}\sigma_{j}+a_{i}a_{j}-\Gamma_{ij}^{k}\sigma_{k}-\Gamma_{ij}^{k}a_{k}\Big]\nonumber\\
&+\sum_{i,j}g^{ij}\Big(c(\partial_{i})A+Ac(\partial_{i})\Big)\partial_{j}
+\sum_{i,j}g^{ij}\Big[c(\partial_{i})\partial_{j}(A)+c(\partial_{i})\sigma_{j}A+c(\partial_{i})a_{j}A\nonumber\\
&+Ac(\partial_{i})\sigma_{j}+Ac(\partial_{i})a_{j}\Big]-\frac{1}{8}\sum_{i,j,k,l}R_{ijkl}\overline{c}(\widetilde{e_i})\overline{c}(\widetilde{e_j})c(\widetilde{e_k})c(\widetilde{e_l})+\frac{1}{4}s+A^2\bigg\}.\nonumber
\end{align}
Then, we obtain
\begin{lem} The following identities hold:
\begin{align}
\sigma_2({{D}_{A}}^3)
&=\sum_{i,j,l}c(dx_{l})\partial_{l}(g^{ij})\xi_{i}\xi_{j} +c(\xi)(4\sigma^k+4a^k-2\Gamma^k)\xi_{k}+2[|\xi|^2A-c(\xi)Ac(\xi)]\\
&+\frac{1}{4}|\xi|^2\sum_{s,t,l}\omega_{s,t}(\widetilde{e_l})[c(\widetilde{e_l})\overline{c}(\widetilde{e_s})\overline{c}(\widetilde{e_t})-c(\widetilde{e_l})c(\widetilde{e_s})c(\widetilde{e_t})]+|\xi|^2A;\nonumber\\
\sigma_{3}({{D}_{A}}^3)&=ic(\xi)|\xi|^{2}.\nonumber
\end{align}
\end{lem}
Write
\begin{eqnarray}
\sigma({{D}_{A}}^3)&=&p_3+p_2+p_1+p_0;
~\sigma({{D}_{A}}^{-3})=\sum^{\infty}_{j=3}q_{-j}.
\end{eqnarray}
By the composition formula of pseudodifferential operators, we have
\begin{align}
1=\sigma({{D}_{A}}^3\circ {{D}_{A}}^{-3})&=
\sum_{\alpha}\frac{1}{\alpha!}\partial^{\alpha}_{\xi}
[\sigma({{D}_{A}}^3)]{{D}}^{\alpha}_{x}
[\sigma({{D}_{A}}^{-3})] \\
&=(p_3+p_2+p_1+p_0)(q_{-3}+q_{-4}+q_{-5}+\cdots) \nonumber\\
&+\sum_j(\partial_{\xi_j}p_3+\partial_{\xi_j}p_2+\partial_{\xi_j}p_1+\partial_{\xi_j}p_0)
(D_{x_j}q_{-3}+D_{x_j}q_{-4}+D_{x_j}q_{-5}+\cdots) \nonumber\\
&=p_3q_{-3}+(p_3q_{-4}+p_2q_{-3}+\sum_j\partial_{\xi_j}p_3D_{x_j}q_{-3})+\cdots,\nonumber
\end{align}
by (4.52), we have

\begin{equation}
q_{-3}=p_3^{-1};~q_{-4}=-p_3^{-1}[p_2p_3^{-1}+\sum_j\partial_{\xi_j}p_3D_{x_j}(p_3^{-1})].
\end{equation}
By (4.49)-(4.53), we have some symbols of operators.
\begin{lem} The following identities hold:
\begin{align}
\sigma_{-3}({{D}_{A}}^{-3})&=\frac{ic(\xi)}{|\xi|^{4}};\\
\sigma_{-4}({{D}_{A}}^{-3})&=
\frac{c(\xi)\sigma_{2}({{D}_{A}}^{3})c(\xi)}{|\xi|^8}
+\frac{c(\xi)}{|\xi|^8}\Big(|\xi|^2c(dx_n)\partial_{x_n}c(\xi')
-2h'(0)c(dx_n)c(\xi)\nonumber\\
&+2\xi_{n}c(\xi)\partial_{x_n}c(\xi')+4\xi_{n}h'(0)\Big).\nonumber
\end{align}
\end{lem}

When $n=6$, then ${\rm tr}_{\wedge^*T^*M}[\texttt{id}]=64$, where ${\rm tr}$ as shorthand of ${\rm trace}$.
Since the sum is taken over $r+\ell-k-j-|\alpha|-1=-6, \ r\leq-1, \ell\leq -3$, then we have the $\int_{{\partial}{M}}\overline{\Phi}$
is the sum of the following five cases:

~\\
\noindent  {\bf case (a)~(I)}~$r=-1, l=-3, j=k=0, |\alpha|=1$.\\

\noindent By (4.47), we get
 \begin{equation}
\overline{\Phi}_1=-\int_{|\xi'|=1}\int^{+\infty}_{-\infty}\sum_{|\alpha|=1}{\rm trace}
[\partial^{\alpha}_{\xi'}\pi^{+}_{\xi_{n}}\sigma_{-1}({{D}_{A}}^{-1})\times\partial^{\alpha}_{x'}\partial_{\xi_{n}}\sigma_{-3}({{D}_{A}}^{-3})](x_0)d\xi_n\sigma(\xi')dx'.
\end{equation}

\noindent  {\bf case (a)~(II)}~$r=-1, l=-3, |\alpha|=k=0, j=1$.\\

\noindent We notice that
  \begin{equation}
\overline{\Phi}_2=-\frac{1}{2}\int_{|\xi'|=1}\int^{+\infty}_{-\infty} {\rm trace}[\partial_{x_{n}}\pi^{+}_{\xi_{n}}\sigma_{-1}({{D}_{A}}^{-1})\times\partial^{2}_{\xi_{n}}\sigma_{-3}({{D}_{A}}^{-3})](x_0)d\xi_n\sigma(\xi')dx'.
\end{equation}

\noindent  {\bf case (a)~(III)}~$r=-1,l=-3,|\alpha|=j=0,k=1$.\\

\noindent By (4.47), we compute that
 \begin{equation}
\overline{\Phi}_3=-\frac{1}{2}\int_{|\xi'|=1}\int^{+\infty}_{-\infty}{\rm trace}[\partial_{\xi_{n}}\pi^{+}_{\xi_{n}}\sigma_{-1}({{D}_{A}}^{-1})
      \times\partial_{\xi_{n}}\partial_{x_{n}}\sigma_{-3}({{D}_{A}}^{-3})](x_0)d\xi_n\sigma(\xi')dx'.
\end{equation}\\
By Lemma 4.2 and Lemma 4.5, we have $\sigma_{-3}(({D}^*_{A}{D}_{A}{D}^*_{A})^{-1})=\sigma_{-3}({{D}_{A}}^{-3})$, then we obtain
\begin{align}
\overline{\Phi}_1+\overline{\Phi}_2+\overline{\Phi}_3=\frac{10}{2}\pi h'(0)\Omega_{4}dx',
\end{align}
 where ${\rm \Omega_{4}}$ is the canonical volume of $S^{4}.$\\

\noindent  {\bf case (b)}~$r=-1,l=-4,|\alpha|=j=k=0$.\\

\noindent By (4.47), we notice that
 \begin{align}
\overline{\Phi}_4&=-i\int_{|\xi'|=1}\int^{+\infty}_{-\infty}{\rm trace}[\pi^{+}_{\xi_{n}}\sigma_{-1}({{D}_{A}}^{-1})
      \times\partial_{\xi_{n}}\sigma_{-4}({{D}_{A}}^{-3})](x_0)d\xi_n\sigma(\xi')dx'\\
&=i\int_{|\xi'|=1}\int^{+\infty}_{-\infty}{\rm trace} [\partial_{\xi_n}\pi^+_{\xi_n}\sigma_{-1}({{D}_{A}}^{-1})\times
\sigma_{-4}({{D}_{A}}^{-3})](x_0)d\xi_n\sigma(\xi')dx'.\nonumber
\end{align}

In the normal coordinate, $g^{ij}(x_{0})=\delta^{j}_{i}$ and $\partial_{x_{j}}(g^{\alpha\beta})(x_{0})=0$, if $j<n$; $\partial_{x_{j}}(g^{\alpha\beta})(x_{0})=h'(0)\delta^{\alpha}_{\beta}$, if $j=n$.
So by  \cite{Wa3}, when $k<n$, we have $\Gamma^{n}(x_{0})=\frac{5}{2}h'(0)$, $\Gamma^{k}(x_{0})=0$, $\delta^{n}(x_{0})=0$ and $\delta^{k}=\frac{1}{4}h'(0)c(\widetilde{e_k})c(\widetilde{e_n})$. Then, we obtain

\begin{align}
\sigma_{-4}({{D}_{A}}^{-3})(x_{0})
&=\frac{1}{|\xi|^8}c(\xi)\Big(h'(0)c(\xi)\sum_{k<n}\xi_kc(\widetilde{e_k})c(\widetilde{e_n})-h'(0)c(\xi)\sum_{k<n}\xi_k\overline{c}(\widetilde{e_k})\overline{c}(\widetilde{e_n})\\
&-5h'(0)\xi_nc(\xi)+2[|\xi|^2A-c(\xi)Ac(\xi)]+\frac{5}{4}|\xi|^2h'(0)c(\widetilde{e_i})\overline{c}(dx_n)\overline{c}(\widetilde{e_i})\nonumber\\
&-\frac{1}{4}|\xi|^2h'(0)c(dx_n)+|\xi|^2A\Big)c(\xi)+\frac{c(\xi)}{|\xi|^8}\Big(|\xi|^2c(dx_n)\partial_{x_n}c(\xi')\nonumber\\
&-2h'(0)c(dx_n)c(\xi)+2\xi_{n}c(\xi)\partial_{x_n}c(\xi')+4\xi_{n}h'(0)\Big).\nonumber
\end{align}
By (4.16) and (4.60), we have
\begin{align}
&{\rm tr} [\partial_{\xi_n}\pi^+_{\xi_n}\sigma_{-1}({{D}_{A}}^{-1})\times
\sigma_{-4}({{D}_{A}}^{-3}) ](x_0)|_{|\xi'|=1} \\
&=-\frac{64h'(0)}{2(\xi_{n}-i)^{2}(1+\xi_{n}^{2})^{4}}\big(-\frac{15}{4}i+16\xi_{n}+\frac{19}{2}i\xi_{n}^{2}+8\xi_{n}^{3}+\frac{21}{4}i\xi_{n}^{4}\big)\nonumber\\
&+\frac{5+6i\xi_{n}-\xi_{n}^{2}}{2(\xi_{n}-i)^{2}(1+\xi_{n}^{2})^{3}}{\rm tr}\big[Ac(\xi')\big]+\frac{-i+6\xi_{n}+5i\xi_{n}^{2}}{2(\xi_{n}-i)^{2}(1+\xi_{n}^{2})^{3}}{\rm tr}\big[Ac(dx_n)\big].\nonumber
\end{align}
By applying the formula shown in (4.59), we can calculate
\begin{align}
\overline{\Phi}_4&=
 i\int_{|\xi'|=1}\int^{+\infty}_{-\infty}\frac{-64h'(0)(-\frac{15}{4}i+16\xi_{n}+\frac{19}{2}i\xi_{n}^{2}+8\xi_{n}^{3}+\frac{21}{4}i\xi_{n}^{4})}{2(\xi_{n}-i)^{2}(1+\xi_{n}^{2})^{4}}d\xi_n\sigma(\xi')dx'\\
 &+i\int_{|\xi'|=1}\int^{+\infty}_{-\infty}\frac{5+6i\xi_{n}-\xi_{n}^{2}}{2(\xi_{n}-i)^{2}(1+\xi_{n}^{2})^{3}}{\rm tr}\big[Ac(\xi')\big]d\xi_n\sigma(\xi')dx'\nonumber\\
 &+i\int_{|\xi'|=1}\int^{+\infty}_{-\infty}\frac{-i+6\xi_{n}+5i\xi_{n}^{2}}{2(\xi_{n}-i)^{2}(1+\xi_{n}^{2})^{3}}{\rm tr}\big[Ac(dx_n)\big]d\xi_n\sigma(\xi')dx'\nonumber\\
&=\frac{111}{2}\pi h'(0)\Omega_4dx'-\frac{5}{16}\pi\Omega_4{\rm tr}\big[Ac(dx_n)\big]dx'.\nonumber
\end{align}

\noindent {\bf  case (c)}~$r=-2,l=-3,|\alpha|=j=k=0$.\\

\noindent We calculate
\begin{equation}
\overline{\Phi}_5=-i\int_{|\xi'|=1}\int^{+\infty}_{-\infty}{\rm trace} [\pi^{+}_{\xi_{n}}\sigma_{-2}({{D}_{A}}^{-1})
      \times\partial_{\xi_{n}}\sigma_{-3}({{D}_{A}}^{-3})](x_0)d\xi_n\sigma(\xi')dx'.
\end{equation}\\
By calculation, we have

\begin{align}
\overline{\Phi}_5=-\frac{105}{4}\pi h'(0)\Omega_4dx'+\frac{161}{512}\pi\Omega_4{\rm tr}[Ac(dx_n)]dx'.
\end{align}

Now $\overline{\Phi}$ is the sum of the cases (a), (b) and (c), hence
\begin{equation}
\overline{\Phi}=\frac{137}{4}\pi h'(0)\Omega_4dx'+\frac{1}{512}\pi\Omega_4{\rm tr}\big[Ac(dx_n)\big]dx'.
\end{equation}

\begin{thm}
Let $M$ be a $6$-dimensional
compact oriented manifold with the boundary $\partial M$ and the metric
$g^M$ as above, ${D}_{A}$ be the perturbation of the de Rham Hodge operator on $\widetilde{M}$, then
\begin{align}
&\widetilde{{\rm Wres}}[\pi^+{{D}_{A}}^{-1}\circ\pi^+({{D}_{A}}^{-3})]=128\pi^3\int_{M}{\rm tr}\bigg(-\frac{1}{12}s+2A^2-\frac{1}{2}\sum_{j=1}^{n}Ac(\widetilde{e_j})Ac(\widetilde{e_j})\bigg)d{\rm Vol_{M}}\nonumber\\
&+\int_{{\partial}{M}}\frac{137}{4}\pi h'(0)\Omega_4d{\rm Vol_{\partial M}}+\int_{{\partial}{M}}\frac{1}{512}\pi\Omega_4{\rm tr}\big[Ac(dx_n)\big]d{\rm Vol_{\partial M}}\nonumber.
\end{align}
where $s$ is the scalar curvature.
\end{thm}

When $A=c(X),$  we can directly state the subsequent corollary:
\begin{cor}
Let $M$ be a $6$-dimensional oriented
compact manifolds with the boundary $\partial M$ and the metric
$g^M$ as above, and let $A=c(X),$ then
\begin{align}
&\widetilde{{\rm Wres}}[\pi^+{{D}_{A}}^{-1}\circ\pi^+({{D}_{A}}^{-3})]=8192\pi^3\int_{M}\bigg(-\frac{1}{12}s\bigg)d{\rm Vol_{M}}\\
&+\int_{\partial M}\frac{137}{4}\pi h'(0)\Omega_4d{\rm Vol_{\partial M}}-\int_{\partial M}\frac{1}{8}\pi\Omega_4g(\partial{x_n}, X)d{\rm Vol_{\partial M}},\nonumber
\end{align}
where $s$ is the scalar curvature.
\end{cor}

When $A=\overline{c}(X),$ we can get the following corollary:
\begin{cor}
Let $M$ be a $6$-dimensional oriented
compact manifolds with the boundary $\partial M$ and the metric
$g^M$ as above, and let $A=\overline{c}(X),$ then
\begin{align}
\widetilde{{\rm Wres}}[\pi^+{{D}_{A}}^{-1}\circ\pi^+({{D}_{A}}^{-3})]&=8192\pi^3\int_{M}\bigg(-\frac{1}{12}s-|X|^2\bigg)d{\rm Vol_{M}}+\int_{\partial M}\frac{137}{4}\pi h'(0)\Omega_4d{\rm Vol_{\partial M}},
\end{align}
where $s$ is the scalar curvature.
\end{cor}

When $A=c(X)c(Y),$ we compute $\widetilde{{\rm Wres}}[\pi^+{{D}_{A}}^{-1}\circ\pi^+({{D}_{A}}^{-3})].$
\begin{cor}
Let $M$ be a $6$-dimensional oriented
compact manifolds with the boundary $\partial M$ and the metric
$g^M$ as above, and let $A=c(X)c(Y),$ then
\begin{align}
&\widetilde{{\rm Wres}}[\pi^+{{D}_{A}}^{-1}\circ\pi^+({{D}_{A}}^{-3})]=8192\pi^3\int_{M}\bigg(-\frac{1}{12}s+|X|^2|Y|^2+4g(X,Y)^2\bigg)d{\rm Vol_{M}}\\
&+\int_{\partial M}\frac{137}{4}\pi h'(0)\Omega_4d{\rm Vol_{\partial M}},\nonumber
\end{align}
where $s$ is the scalar curvature.
\end{cor}

When $A=c(X)\overline{c}(Y),$ we notice that:
\begin{cor}
Let $M$ be a $6$-dimensional oriented
compact manifolds with the boundary $\partial M$ and the metric
$g^M$ as above, and let $A=c(X)\overline{c}(Y),$ then
\begin{align}
&\widetilde{{\rm Wres}}[\pi^+{{D}_{A}}^{-1}\circ\pi^+({{D}_{A}}^{-3})]=8192\pi^3\int_{M}\bigg(-\frac{1}{12}s\bigg)d{\rm Vol_{M}}+\int_{\partial M}\frac{137}{4}\pi h'(0)\Omega_4d{\rm Vol_{\partial M}},
\end{align}
where $s$ is the scalar curvature.
\end{cor}

When $A=\overline{c}(X)\overline{c}(Y),$ similar to (4.68), we can get the following corollary:
\begin{cor}
Let $M$ be a $6$-dimensional oriented
compact manifolds with the boundary $\partial M$ and the metric
$g^M$ as above, and let $A=\overline{c}(X)\overline{c}(Y),$ then
\begin{align}
&\widetilde{{\rm Wres}}[\pi^+{{D}_{A}}^{-1}\circ\pi^+({{D}_{A}}^{-3})]=8192\pi^3\int_{M}\bigg(-\frac{1}{12}s-|X|^2|Y|^2+6g(X,Y)^2\bigg)d{\rm Vol_{M}}\\
&+\int_{\partial M}\frac{137}{4}\pi h'(0)\Omega_4d{\rm Vol_{\partial M}},\nonumber
\end{align}
where $s$ is the scalar curvature.
\end{cor}

When $A=c(X)c(Y)c(Z),$ we compute that:
\begin{cor}
Let $M$ be a $6$-dimensional oriented
compact manifolds with the boundary $\partial M$ and the metric
$g^M$ as above, and let $A=c(X)c(Y)c(Z),$ then
\begin{align}
&\widetilde{{\rm Wres}}[\pi^+{{D}_{A}}^{-1}\circ\pi^+({{D}_{A}}^{-3})]=8192\pi^3\int_{M}\bigg(-\frac{1}{12}s-2|X|^2|Y|^2|Z|^2+2|X|^2g(Y,Z)^2\\
&+2|Y|^2g(X,Z)^2+2|Z|^2g(X,Y)^2-4g(X,Y)g(X,Z)g(Y,Z)\bigg)d{\rm Vol_{M}}\nonumber\\
&+\int_{\partial M}\frac{137}{4}\pi h'(0)\Omega_4d{\rm Vol_{\partial M}}+\int_{\partial M}\frac{1}{8}\pi\Omega_4g(\partial{x_n}, X)g(Y,Z)d{\rm Vol_{\partial M}}\nonumber\\
&-\int_{\partial M} \frac{1}{8}\pi\Omega_4g(\partial{x_n}, Y)g(X,Z)d{\rm Vol_{\partial M}}+\int_{\partial M} \frac{1}{8}\pi\Omega_4 g(\partial{x_n}, Z)g(X,Y)d{\rm Vol_{\partial M}},\nonumber
\end{align}
where $s$ is the scalar curvature.
\end{cor}

When $A=\overline{c}(X)c(Y)c(Z),$ we can conclude the following facts:
\begin{cor}
Let $M$ be a $6$-dimensional oriented
compact manifolds with the boundary $\partial M$ and the metric
$g^M$ as above, and let $A=\overline{c}(X)c(Y)c(Z),$ then
\begin{align}
&\widetilde{{\rm Wres}}[\pi^+{{D}_{A}}^{-1}\circ\pi^+({{D}_{A}}^{-3})]=8192\pi^3\int_{M}\bigg(-\frac{1}{12}s+3|X|^2|Y|^2|Z|^2-4|X|^2g(Y,Z)^2\bigg)d{\rm Vol_{M}}\\
&+\int_{\partial M}\frac{137}{4}\pi h'(0)\Omega_4d{\rm Vol_{\partial M}},\nonumber
\end{align}
where $s$ is the scalar curvature.
\end{cor}

Computations show that:
\begin{cor}
Let $M$ be a $6$-dimensional oriented
compact manifolds with the boundary $\partial M$ and the metric
$g^M$ as above, and let $A=\overline{c}(X)\overline{c}(Y)c(Z),$ then
\begin{align}
&\widetilde{{\rm Wres}}[\pi^+{{D}_{A}}^{-1}\circ\pi^+({{D}_{A}}^{-3})]=8192\pi^3\int_{M}\bigg(-\frac{1}{12}s-4|X|^2|Y|^2|Z|^2+4|Z|^2g(X,Y)^2\bigg)d{\rm Vol_{M}}\nonumber\\
&+\int_{\partial M}\frac{137}{4}\pi h'(0)\Omega_4d{\rm Vol_{\partial M}}-\int_{\partial M}\frac{1}{8}\pi\Omega_4g(\partial{x_n}, Z)g(X, Y)d{\rm Vol_{\partial M}},\nonumber
\end{align}
where $s$ is the scalar curvature.
\end{cor}

When $A=\overline{c}(X)\overline{c}(Y)\overline{c}(Z),$ we can get the following corollary:
\begin{cor}
Let $M$ be a $6$-dimensional oriented
compact manifolds with the boundary $\partial M$ and the metric
$g^M$ as above, and let $A=\overline{c}(X)\overline{c}(Y)\overline{c}(Z),$ then
\begin{align}
&\widetilde{{\rm Wres}}[\pi^+{{D}_{A}}^{-1}\circ\pi^+({{D}_{A}}^{-3})]=8192\pi^3\int_{M}\bigg(-\frac{1}{12}s+5|X|^2|Y|^2|Z|^2-6|X|^2g(Y,Z)^2\\
&-6|Y|^2g(X,Z)^2-6|Z|^2g(X,Y)^2+12g(X,Y)g(X,Z)g(Y,Z)\bigg)d{\rm Vol_{M}}\nonumber\\
&+\int_{\partial M}\frac{137}{4}\pi h'(0)\Omega_4d{\rm Vol_{\partial M}},\nonumber
\end{align}
where $s$ is the scalar curvature.
\end{cor}

\vskip 1 true cm

\section{Acknowledgements}

The author was supported in part by  NSFC No.11771070. The author thanks the referee for his (or her) careful reading and helpful comments.

\vskip 1 true cm


\bigskip
\bigskip

\noindent {\footnotesize {\it S. Liu} \\
{School of Mathematics and Statistics, Northeast Normal University, Changchun 130024, China}\\
{Email: liusy719@nenu.edu.cn}

\noindent {\footnotesize {\it T. Wu} \\
{School of Mathematics and Statistics, Northeast Normal University, Changchun 130024, China}\\
{Email: wut977@nenu.edu.cn}

\noindent {\footnotesize {\it Y. Wang} \\
{School of Mathematics and Statistics, Northeast Normal University, Changchun 130024, China}\\
{Email: wangy581@nenu.edu.cn}

\end{document}